\documentclass{amsart}
\usepackage{psfig,latexsym,amssymb}
\def\proof{\medskip\noindent{\sc Proof. }}
\def\dps{\Delta(\Pi_n)/S_n}
\def\dwreath{\Delta(B_{2n})/S_2\wr S_n }

\def\EOP{\hfill$\Box$}

\def\natnum{\hbox{\rm I\kern-.17em N}}
\def\integ{\hbox{\rm Z\kern-.3em Z}}
\def\reals{{\mathbb R}}

\newtheorem{thm}{Theorem}[section]
\newtheorem{lem}{Lemma}[section]
\newtheorem{prop}{Proposition}[section]
\newtheorem{defn}{Definition}[section]
\newtheorem{cor}{Corollary}[section]
\newtheorem{qn}{Question}[section]
\newtheorem{examp}{Example}[section]

\newtheorem{rk}{Remark}[section]
\newtheorem{conj}{Conjecture}[section]
\newtheorem{partcon}{Partitioning Construction}[section]
\begin{document}

\title{Lexicographic Shellability for Balanced Complexes}
\author{Patricia Hersh}
\address{Department of Mathematics, Box 354350,
University of Washington, Seattle, WA 98195}
\curraddr{
        Department of Mathematics,
	525 East University Ave.,
        University of Michigan,
        Ann Arbor, Michigan 48109-1109}
\email{plhersh@umich.edu}
%\date{February 18th, 2000}

\begin{abstract}
We introduce a notion of lexicographic shellability for pure, balanced 
boolean cell complexes, modelled after the $CL$-shellability criterion 
of Bj\"orner and Wachs for posets [BW2] 
and its generalization by Kozlov [Ko2] called
$CC$-shellability.    
We give a lexicographic shelling for the quotient of the order complex
of a Boolean algebra
of rank $2n$ by the action of the wreath product $S_2\wr S_n$ of symmetric
groups, and we provide a partitioning for the quotient complex
$\Delta (\Pi_n )/S_n $.  

Stanley asked for a description of the symmetric group representation
$\beta_S $ on the homology of the rank-selected partition lattice
$\Pi_n^S $
in [St2], and in particular he asked when the multiplicity $b_S(n)$ of the
trivial representation in $\beta_S$ is $0$.  
One consequence of the partitioning for $\dps $ is 
a (fairly complicated) combinatorial interpretation for 
$b_S(n) $; another is a simple proof of Hanlon's result [Ha] that 
$b_{1,\dots ,i}(n)=0$.  Using a result of Garsia and Stanton from [GS], 
we deduce from our shelling for 
$\Delta (B_{2n})/S_2 \wr S_n$ that the ring of invariants
$k[x_1,\dots ,x_{2n}]^{S_2\wr S_n}$ is Cohen-Macaulay over any field $k$.
\end{abstract}

\maketitle

\section{Introduction}\label{introsec}

Let $B_n$ denote the Boolean algebra of subsets of $\{ 1,\cdots ,n\} $ 
ordered by inclusion and let $\Pi_n $ be the lattice of unordered partitions 
of $\{ 1,\cdots ,n\} $ ordered by refinement.  The natural symmetric 
group action on $\{ 1,\cdots ,n \} $ induces an action on each of these 
posets.  Likewise, the wreath product $S_2 \wr S_n $ acts on the 
elements of the Booolean algebra $B_{2n}$.  Any rank-preserving 
group action on a finite, ranked poset $P$ with 
minimal and maximal elements $\hat{0} $ and $\hat{1} $ induces an 
action on the order complex $\Delta(P)$, that is, on the simplicial complex 
consisting of an $(i-1)$-face for each $i$-chain $\hat{0} < u_1 <
\cdots < u_i < \hat{1} $ in $P$; the group action on poset elements induces
an action on chains.  
This gives rise to a quotient cell complex, denoted $\Delta(P)/G$, which is
comprised of the $G$-orbits of order complex faces.  

Note that $\Delta(P)/G $ need not
coincide with the order complex of the quotient poset $P/G$ because
there may be covering relations $u<v$ and $u'<v'$ in $P$ belonging to 
distinct orbits despite having $u'=gu$ and $v'=g'v$ for some $g,g'\in G$.
Babson and Kozlov give
conditions under which $\Delta(P)/G = \Delta(P/G)$ in [BK].  Equality 
does not hold for $P=\Pi_n, G=S_n$ and for $P=B_{kn}, G=S_k\wr S_n$, so
the quotient complexes $\Delta (\Pi_n)/S_n$ and $\Delta (B_{2n})/S_2\wr S_n$
that we will consider are not simplicial complexes and in particular
are not order complexes of posets.

We shall give a lexicographic shelling for the quotient complex 
$\Delta(B_{2n})/S_2\wr S_n$ and a partitioning for
$\Delta (\Pi_n )/S_n$, using a generalized notion of chain-labelling for
balanced complexes.  By way of comparison, Ziegler showed 
in [Zi] that the quotient poset $\Delta(\Pi_n /S_n ) $ is not 
Cohen-Macaulay for $n\ge 19$.  We will verify in 
Section \ref{wreathsec} that $\Delta(B_6)/S_3\wr S_2$ is not shellable.  
It was shown in [GS], [Re], [St2] that shellings and partitionings
for quotient complexes $\Delta(P)/G$ yield
information about sub-rings of invariant polynomials 
and about $G$-representations on the 
homology of $\Delta(P)$.  Our approach is
to extend poset lexicographic shelling methods to a general enough class
of complexes to include quotient complexes $\Delta(P)/G$, taking advantage
of properties quotient complexes have in common with posets.  

To this end, we introduce a notion of 
$CC$-labelling for pure, balanced boolean cell 
complexes in Section \ref{lexsec} and confirm that it induces a 
lexicographic shelling.  We require three conditions for a 
chain-labelling to be a $CC$-labelling.  The first of these is a direct
translation of the poset requirement that each interval must have a unique
(topologically) increasing chain and
that this be lexicographically smallest on the interval.  The
second condition, which we call the ``crossing condition'', is automatic
for posets, and we verify that it also holds for all quotient 
complexes.  The
third condition, the ``multiple-face-overlap'' condition, is vacuous for 
simplicial complexes.  In the quotient complex $\Delta(B_{2n})/S_2\wr S_n$
this technical condition follows readily from the increasing chain
condition because of the nature of the ascents that occur in 
the labelling for $\dwreath $.  A virtually identical argument isolates
exactly where our lexicographic order on facets in $\dps $ fails to 
satisfy the multiple-face-overlap condition.  It also helps us find a face
whose link is the real projective plane, implying $\dps $ is not
Cohen-Macaulay over $\integ /2\integ $, and hence is not shellable.

Section \ref{wreathsec} provides a
$CC$-labelling for $\dwreath $ while Section \ref{partsec} gives a 
chain-labelling for $\dps $ that satisfies the increasing chain condition.  
Section ~\ref{partsec} constructs from this labelling 
a partitioning for $\dps $.  In Section \ref{trivsec},
we express the multiplicity $b_S(n)$ of the trivial
representation in the $S_n$-representation on the homology of 
the rank-selected complex $\Delta(\Pi_n^S) $ in terms of the flag
$h$-vector for the quotient complex $\Delta(\Pi_n)/S_n $.  
This gives a combinatorial interpretation
for $b_S(n)$, since $h_S (\dps )$ is the number of facets contributing
minimal new faces of support $S$.  In any lexicographic shelling (or a
partitioning which uses descents in a similar spirit), this 
support $S$ is the set of topological descents in each facet.  One 
may use this to show that $b_S(n)>0$ for a particular $S$ by showing that 
the descent set $S$ arises in some lexicographic shelling (or 
partitioning) step (cf. [HH] for more results of this nature).  
In Section 7, an analysis of which $S$ arise as minimal new faces in a 
partitioning yields a simple proof of 
Hanlon's result that $b_{1,\dots ,i}(n)=0$.   Finally, Section 
\ref{invtsec} gives an application to subrings of invariant polynomials 
by applying a result of Garsia and Stanton in [GS].

The remainder of our introduction will review the notion of boolean cell
complex, lexicographic shellability for posets and related
terminology.  Boolean cell complexes were introduced by Bj\"orner in [Bj] 
and by Garsia and Stanton in [GS]. Stanley studied their face posets,
namely simplicial posets, in [St3].  He defined the face ring 
(also called Stanley-Reisner ring) of
a boolean cell complex and then showed that a boolean cell 
complex is Cohen-Macaulay if and only if its face ring is
Cohen-Macaulay.  Duval studied free resolutions of face rings
of boolean cell complexes in [Du].  Reiner developed a theory of
$P$-partitions for Coxeter groups in order to shell (or in some
cases partition) quotients of Coxeter complexes in [Re].  Our 
interest is in lexicographic shelling for balanced, boolean cell complexes.

\begin{defn}[Bj\"orner, Garsia-Stanton]
A regular cell complex is {\bf boolean} if every interval
in its face poset is a boolean algebra.  
\end{defn}

A boolean cell complex is much like a simplicial complex,
except that more than one face may have the same set of vertices.
We note that boolean cell complexes may alternatively
be defined quite naturally in terms of simplicial sets.
We find it convenient to refer to $i$-cells as $i$-faces,
0-cells as vertices, and 
in the same vein to call cells of top dimension facets.
A boolean cell complex is {\bf pure} of dimension $n$
if all the maximal cells have dimension $n$,
and then it 
is {\bf balanced } if the vertices may be colored with 
$n+1$ colors so that no two vertices in a face have the 
same color.  We refer to the set of colors for the vertices in a face 
as the {\bf support} of the face and say a face has {\bf disconnected}
support if the support includes $i,k$ for $i<k$ and does not include
some $j$ for $i<j<k$.  

Bj\"orner established the following notion of
shellability (phrased slightly differently) for 
boolean cell complexes in [Bj1]. 
\begin{defn}[Bj\"orner]
A boolean cell complex is {\bf shellable} if the 
facets may be ordered $F_1,\dots ,F_k$ so that 
$F_j \cap (\cup_{i=1}^{j-1} F_i )$ is pure of codimension one for each
$1<j\le k$. 
\end{defn}
Recall that the order complex of a finite, ranked poset is a 
balanced simplicial
complex in which vertices are colored by poset rank.  This will 
allow us to translate
poset notions of lexicographic shellability to conditions on 
the order complex which may then be extended to give shelling criteria for
more general balanced complexes.

One may conclude from the 
existence of a shelling that a boolean complex has the homotopy type of a 
wedge of spheres of top dimension and is 
Cohen-Macaulay (cf. [Bj1]).  The arguments are similar (though slightly more
subtle) to those for simplicial complexes.

\begin{prop}[Bj\"orner]\label{goofsphere}
If a pure boolean cell complex is shellable, then it has the 
homotopy type of a wedge of spheres of top dimension.
\end{prop}

In a shelling, each facet either attaches along its entire boundary, 
``closing off'' a sphere, or its overlap with the union of
earlier facets is a simplicial complex with a cone point, implying 
the homotopy type is unchanged by the facet's insertion.

\begin{prop}[Bj\"orner]\label{goofcm}
If a pure boolean cell complex is shellable, then it is Cohen-Macaulay.
\end{prop}

This may be shown by shelling its barycentric subdivision (cf. [Bj2, p. 173]),
which is a simplicial complex, then invoking Munkres' classical result 
[Mu2, p.117,121-123] that Cohen-Macaulayness does 
not depend on choice of triangulation.

\begin{defn}[Bj\"orner]
An integer labelling of the covering relations in a 
finite poset with $\hat{0} $ and $\hat{1}$
is an {\bf EL-labelling} if it has the following two
properties, which together constitute the {\bf increasing chain condition}.
\begin{enumerate}
\item
Every interval has a unique saturated chain with (weakly) increasing
edge labels.
\item
The increasing chain is the lexicographically smallest chain of 
labels on an interval.
\end{enumerate}
\end{defn}

A {\bf chain-labelling} is a labelling of poset covering 
relations such that the label
assigned to a covering relation $u\prec v$ may depend on the 
choice of root, namely on the saturated chain $\hat{0}\prec u_1\prec\cdots
\prec u_k=u$ as well as on $u$ and $v$.  Recall from [BW1], [BW2] that
a {\bf CL-labelling} is any chain-labelling satisfying the increasing-chain 
condition.  An edge-labelling or chain-labelling induces a partial order
on facets by lexicographically ordering the sequences of labels assigned to
saturated chains.  It is shown in [Bj1] (resp. [BW2]) that any total order
extension of the lexicographic order given by an EL-labelling (resp. 
CL-labelling) is a shelling order on facets.

Kozlov generalized poset EL-shellability (resp. CL-shellability) in [Ko2]
to a criterion he called {\bf EC-shellability} (resp. {\bf CC-shellability}) by
relaxing the requirement that every poset interval must 
have a unique increasing chain.  In effect, he instead
requires each interval to have a unique saturated chain
that behaves topologically like an increasing chain with respect to the 
chosen lexicographic order.
We rediscovered EC/CC-shellability in the course of joint work with 
Kleinberg (see [HK]); it will be convenient for the fairly involved 
shelling arguments in later sections to use the notation and point of 
view taken in [HK], so now let us review terminology from [HK].

We classify ascents and descents in a poset 
chain-labelling (or edge-labelling) $\lambda $ as follows: 
let us say that a pair of edges $u\prec v$ and $v \prec w$ constitute a 
{\bf topological ascent} if the word consisting of two consecutive
labels $\lambda (u,v) $ and $\lambda(v,w)$ is lexicographically
smallest on the interval from $u$ to $w$ and let us say that the pair of
covering relations $u\prec v$ and $v\prec w$ comprise a {\bf topological 
descent} otherwise.  We may further distinguish between topological 
ascents with increasing or decreasing consecutive labels by calling the 
former {\bf honest ascents} and the latter {\bf swap descents}.  
Similarly, we call 
topological descents with decreasing labels {\bf honest descents} and all 
others {\bf swap ascents}.  In this language, a poset is EC-shellable 
(resp. CC-shellable) if each interval (resp. rooted interval) has a unique
topologically increasing chain (namely a chain consisting entirely of 
topological ascents) and this is the lexicographically smallest chain on
the interval.  It is shown in [Ko2] that these labellings induce 
lexicographic shellings, for the same reason that EL-labellings and 
CL-labellings do.

\section{A lexicographic shelling condition for balanced 
boolean cell complexes}\label{lexsec}

In this section, we extend CC-shellability to pure, balanced 
boolean cell complexes and EL/CL/EC/CC-shellability to
pure, balanced simplicial complexes.  We always choose indices so that
$F_i$ precedes $F_j$ lexicographically for $i<j$.  When the vertices 
$v_1,\dots ,v_t$ in a face $\sigma $ are colored $c_1,\dots ,c_t$, then 
we call the set $S=\{ c_1,\dots ,c_t \}$ the {\bf support} of $\sigma $.  
It will be convenient to represent an arbitrary color set as
$s_0,\dots ,r_1,s_1,\dots ,r_2,s_2,\dots ,r_k,s_k,\dots ,
r_{k+1}$ for some $s_0\ge 1$, $r_{k+1}< n$ and with $s_i-r_i >1$ 
for all $1\le i\le k$.  We use 
$s_i,\dots ,r_{i+1}$ to denote the collection of all
possible colors from $s_i$ to $r_{i+1}$.  Thus, the colors not in $S$ are 
those between $r_i$ and $s_i$ for some $i$, 
and the colors that are smaller than $s_0$ or larger than $r_{k+1}$.

Let us begin by translating the poset lexicographic
shellability condition of Bj\"orner and Wachs [BW1]
to a condition on the order complex so as to make an analogous
condition for pure, balanced boolean cell complexes.  Consider any
finite, graded poset with unique minimal and maximal elements $\hat{0} $ and 
$\hat{1}$.  Notice that the increasing chain condition 
on the Hasse diagram of a poset may be viewed
as a condition on the order complex.  We make the following conventions:
 
\begin{enumerate}
\item
The label assigned to
a poset covering relation from rank $i$ to $i+1$ is placed on the consequent
order complex edge colored $i,i+1$ for $1\le i \le n-1$.
\item
The label of each poset edge involving $\hat{0} $
(resp. $\hat{1}$) is assigned to the corresponding
vertex colored 1 (resp. $n-1$).  
\item
The interval from $u$ to $v$ in the poset is the collection of faces 
colored $rk(u),rk(u)\\
+1,\dots ,rk(v)-1,rk(v)$ 
which include the vertices $u$ and $v$.
\item
Let $i=rk(u)$ and $j=rk(v)$.  Then
each poset saturated chain from $u$ to $v$
translates to a walk on the resulting face colored $i,i+1,\dots ,j-1,j$
in the order complex.  Such a walk along edges colored $i',i'+1$ for 
$1\le i'<j$ passes through the vertex colors sequentially.
\end{enumerate}

The increasing-chain condition on an interval from rank
$i$ to $j$ amounts to
a condition on all the faces in the order complex
consisting of vertices colored $i,\dots ,j$
which include a particular pair of vertices colored $i$ and $j$.  
This requirement 
makes sense for arbitrary pure, balanced simplicial complexes, using the 
balancing to play the role of poset rank.  Any pure, balanced simplicial
complex will be lexicographically shellable if it satisfies this
increasing chain condition along with another requirement which we call
the crossing condition.
For balanced boolean cell complexes, we must define the notion of
interval a little bit more carefully, but again the increasing chain 
condition will be a similar requirement on cells in an interval; the 
increasing chain condition together with the crossing condition and a third
requirement called the multiple-face-overlap condition will imply that a 
pure, balanced boolean cell complex is lexicographically shellable.

Let us generalize the notions of interval and rooted interval to 
balanced complexes, as follows.

\begin{defn}
Let $\tau $ be a face colored $1,\dots ,i,j$ for some $i<j$ in a pure, 
balanced boolean cell complex.  The 
{\bf rooted interval} specified by $\tau $ is the collection of faces
colored $1,\dots ,j$ that contain $\tau $.  
\end{defn}

Notice that a cell complex which is not a simplicial complex might have
several faces comprised of the same vertex set of support $1,\dots ,i,j$.
Each of these faces gives rise to a different interval.  For this reason,
it does not seem wise to allow edge-labellings and un-rooted intervals
when workin with boolean cell complexes that are not simplicial complexes.

\begin{defn}
Let $u,v $ be a pair of vertices colored $i,j$ for $i<j$ in a pure, 
balanced simplicial complex.  Then (unrooted) {\bf interval} specified by
$u$ and $v$ is the collection of faces colored $i,\dots ,j$ which contain
the vertices $u$ and $v$.
\end{defn} 

Next, we adapt the 
definition of topological ascent and descent to balanced boolean cell 
complexes.

\begin{defn}
A facet $F_j$ 
has a {\bf topological descent} at the color $r$ if there is a codimension one 
face in $F_j \cap (\cup_{i<j} F_i )$ omitting only the vertex
colored $r$.  Otherwise, it has a {\bf topological ascent} at rank $r$. 
\end{defn}

In a poset, one may view the saturated chains as 
non-self-intersecting paths from $\hat{0} $ to $\hat{1} $ 
in the Hasse diagram, so then saturated chains
intersect where two of these paths cross each other.  
Notice that when two poset saturated chains cross $c-1$ times in the 
proper part of the poset, one obtains $2^c$ distinct saturated chains
by choosing which of the two saturated chains to follow on each of the 
$c$ segments between consecutive crossing points.  The existence of these
poset chains implies for lexicographic orders that every maximal face in
$F_j\cap (\cup_{i<j} F_i)$ skips a single interval of consecutive ranks.
The crossing condition is designed to test for this behavior 
in arbitrary balanced complexes.  
%\begin{figure}[h]
%\begin{picture}(150,90)(-70,0)
%\psfig{figure=cross.eps,height=3cm,width=1cm}
%\end{picture}
%\caption{Two crossings yielding four saturated chains}
%\label{cross}
%\end{figure}

Notice in the case of boolean cell complexes that are not simplicial
complexes that the crossing condition, given next, does not always ensure 
that maximal faces in $F_j\cap (\cup_{i<j} F_i )$ have support skipping a 
single interval of consecutive colors.  Specifically, it does not apply
to faces $\sigma \in F_j\cap F_i $ of support $1,\dots ,r_1,s_1,\dots ,r_2,
\dots ,s_k,\dots ,r_{k+1}$ such that another face in $F_j\cap F_i$ has 
support $1,\dots ,r+1$.  The multiple-face-overlap condition accounts for
these faces.

\begin{defn}\label{cc-shell}
A balanced boolean cell 
complex is {\bf CC-shellable} if there is a chain-labelling satisfying the 
following three conditions.
\begin{enumerate}
\item
{\bf Increasing chain condition.}  
Each rooted edge-interval has a unique extension
with topologically increasing labels and this is the
lexicographically smallest face in the interval (i.e. its lexicographically
earliest extension to a facet is lexicographically smallest among facets that
may be obtained from faces in the interval).
\item
{\bf Crossing condition.}  Let 
$\sigma $ be a face in the intersection of a facet $F_k$ with a 
lexicographically earlier facet $F_j$.  Suppose that (1)
the support of $\sigma $ includes $1,\dots ,r$ while no other face in 
$F_j\cap F_k$ includes support $1,\dots ,r'$ for $r'>r$ and 
(2) the complement of $\sigma $ has disjoint support.   Then there is some 
facet $F_i$ for $i<k$ and some face $\tau \in F_i \cap F_k$ 
such that $\sigma \subsetneq \tau $ and the complement of 
$\tau $ has connected support $r+1,\dots ,s$ for some $s\ge r+1$.
\item
{\bf Multiple-face-overlap condition.}  Suppose the intersection of
two facets $F_i$ and $F_j$ (with $i<j$ ) contains two
faces $\sigma ,\tau $, such that $\sigma $ is maximal in 
$F_j\cap (\cup_{i'=1}^{j-1} F_{i'} )$ and $\tau $ is maximal 
in $F_j\cap F_i$.  
Furthermore, assume 
that $\sigma $ has support including the colors $1,\dots ,r'$ for some
$1,\dots ,r'$ which is not a subset of the support of $\tau $.  Then 
$\tau $ must be contained in a codimension one face $\gamma $ of $F_j$
such that $\gamma\in F_j \cap (\cup_{i'=1}^{j-1} F_i')$.  
Letting $1,\dots ,r,s,\dots ,n$ denote the support of $\tau $, 
it suffices to check this for $s\le r'$.
\end{enumerate}
\end{defn}

The final remark in the multiple-face-overlap condition will be invaluable
to our proofs in later sections and is confirmed within the proof of 
Theorem \ref{shellcrit}.  It allows
us to assume when some $F_i \cap F_j$ includes maximal faces
$\sigma $ and $\tau $ as above that the first covering relation of $F_j$ 
skipped in $\tau $ has larger label than the covering relation of the 
same rank in $F_i$ (for $i<j$).  
To prove such a face $\sigma $ has codimension one,
we may assume (to get a contradiction) that the interval of 
$F_j$ skipped by $\sigma $ consists entirely of topological ascents.

\begin{rk}
{\rm
The above criterion specializes to pure, balanced simplicial 
complexes, in which case the multiple-face-overlap condition is 
vacuously true.  For balanced simplicial complexes, the above criterion is
easily modified to give notions of EL/CL/EC-shellability as follows.  For
EL/CL-shellability, we require increasing chains instead of merely 
topologically increasing chains.  For EL/EC-shellability, we label edges
in a way that does not depend on the root, and intervals are specified by 
pairs of vertices rather than also depending on the entire root.
}
\end{rk}

We call any ordering
of the facets of a balanced complex which is induced by an 
EL/CL/EC/CC-labelling a {\bf lexicographic shelling}.  Let us check
that these do indeed give shellings.

\begin{thm}\label{shellcrit}
If $F_1,\dots ,F_r$ is a lexicographic shelling for a pure,
balanced boolean cell complex $\Delta $ of dimension $n-1$, then
$F_l \cap \left( \cup_{k=1}^{l-1} F_k \right) $ is pure of codimension 
one for each $l$, so $F_1,\dots ,F_r$ is a shelling.
\end{thm}

\proof
Let $H$ be a maximal face in $F_j \cap (\cup_{i<j} F_i)$, so $H\subseteq
F_j\cap F_{i'}$ for some $i'<j$.  It suffices to show that $H$ has codimension
one in $F_j$.  Let $S= \{ s_0,\dots ,r_1,s_1,\dots ,
r_{k+1} \} $ be the support of $H$.  The crossing condition implies that
the complement of $S$ is a single interval of consecutive ranks, except in the
following scenario.
The crossing condition does not apply if 
(1) the support of $H$ includes $1,\dots ,r$ but excludes $r+1$, and 
(2) there is some other face
$H'\subseteq F_j \cap F_{i'}$ has support including $1,\dots ,r+1$.  However,
in this case the multiple-face-overlap condition ensures that $H$ has 
codimension one in $F_j$, as desired.  Hence, we only need to consider $H$ 
of support $S = \{ 1,\dots ,r,s,\dots ,n \}$ or $S=\{ s_0 ,\dots ,n\} $ 
for some $s_0>1 $ or $S=\{ 1,\dots ,r_1 \} $ for some $r_1 < n$.
Let us assume $H$ has support
$1,\dots ,r,s,\dots ,n$, since the other arguments are similar.
The facet $F_{i'}$ must be
strictly smaller in lexicographic order than $F_j$ on the interval skipped 
by $H$, since $F_j\cap F_{i'}$ does not include any faces of support 
$1,\dots ,r+1$.  Thus, the increasing chain condition ensures that $F_j$ 
must have a topological descent on the rooted interval specified by $H$ 
restricted to color set $1,\dots ,r,s$.  Let $F_{i''}$ be the facet in 
which one such topological descent is replaced by a 
topological ascent.  $F_{i''}$ precedes $F_j$ in lexicographic
order, and $F_j\cap F_{i''}$ contains a face of codimension one in 
$F_j$ which contains $H$.  Hence, $H$ must be codimension one in $F_j$ to 
be maximal in $F_j \cap (\cup_{i<j} F_i)$.

Next, we verify the remark in the multiple-face-overlap condition,
recalling that $\tau \subseteq F_i\cap F_j$ for some $F_i$ preceding $F_j$
in lexicographic order. 
Since $\tau $ is assumed to have support 
$1,\dots ,r,s,\dots ,n$, the increasing chain
condition implies that either 
$\tau $ has a topological descent on the interval from color $r$ to $s$, or
$\tau $ is lexicographically smallest on this interval.  In the former case,
we argue as above.  In the latter case, $F_i\cap F_j$ must agree up to color
$s$, namely there must be a face $\sigma\subseteq F_i\cap F_j$ whose support
includes $1,\dots ,s$, just as asserted.
\EOP

\begin{rk}
{\rm
Just as in a poset lexicographic shelling, the topologically
decreasing chains give rise to 
facets attaching along their entire boundaries, and the 
homotopy type of a lexicographically shellable balanced boolean cell
complex is a wedge of spheres of top
dimension, indexed by the topologically decreasing chains.
}
\end{rk}

One could define $CL$-shellability for balanced complexes by 
requiring all of the topological ascents (resp. topological descents) 
in a $CC$-shelling to be actual ascents 
(resp. descents).  This, however, would give a more restrictive notion
of $CL$-shellability than the one very recently introduced 
by Hultman in [Hu].

We conclude this section by confirming that quotient complexes
always satisfy the crossing condition.  Later sections will give a
lexicographic shelling for the quotient complex
$\dwreath $ and use lexicographic shelling ideas in a partitioning for
$\dps $.  Lexicographic shellability for another class of balanced
boolean cell complexes, the nerves of ranked, loop-free small categories,
is discussed in [HK].

\begin{prop}\label{quot}
Quotient complexes $\Delta(P)/G$ satisfy the crossing condition.
\end{prop}

\proof
Suppose that two saturated chain
orbits $F_j, F_k$ (with $j<k$)
share a maximal face $\sigma $ of support $S = \{ s_0,\dots ,r,
s_1,\dots ,r_2,s_2,\dots ,r_{k+1} \} $.  We may assume $s_0=1$, by 
assumption (1) of the crossing condition.
Then there must be poset saturated chains $C_j, C_k$ belonging to the orbits 
$F_j, F_k$, respectively, such that $C_j, C_k$ also share a 
face of support $S$.    
Now consider the saturated chain $C_i$ 
in $P$ which agrees with $C_j$ on ranks $1,\dots ,s_1$ and agrees with
$C_k$ on ranks $s_1,\dots ,n$.  Denote the orbit of $C_i$ by
$F_i$.  We may assume that $F_j$ and $F_k $ do not share a face
colored $\{ 1,\dots ,r+1\} $, by assumption (1) 
of the crossing condition.  Thus, $F_j$ already has a strictly earlier 
labelling than $F_k$ on the interval up to the color $s_1$.  Since 
$F_i$ agrees with $F_j$ on color set $1,\dots ,s_1$, $F_i$ also must 
precede $F_k$ lexicographically.  
Restricting $F_i$ to color set $1,\dots ,r,s_1,\dots ,n$ gives the
desired face $\tau $.
\EOP

\section{A lexicographic shelling for $\Delta(B_{2n})/S_2\wr S_n $}
\label{wreathsec}

\begin{figure}[h]
\begin{picture}(150,62)(-65,0)
\psfig{figure=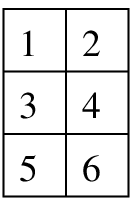,height=2cm,width=1.3cm}
\end{picture}
\caption{Labelled boxes acted upon by $S_2\wr S_3$}
\label{wreathbox}
\end{figure}

An edge-labelling for the lattice $B_{2n}$ of subsets of 
$\{ 1,\dots ,2n \} $ ordered by inclusion
comes from labelling each covering relation
$\{ \sigma_1,\dots ,\sigma_{i-1} \} \subseteq \{ \sigma_1,\dots ,
\sigma_{i-1} ,\sigma_i \}$ with the number $\sigma_i\in \{ 1,\dots ,
2n\} $ being inserted.  This 
labelling assigns the permutation $\sigma_1\sigma_2\cdots
\sigma_{2n}\in S_{2n}$ to the saturated chain
$\emptyset\prec \{ \sigma_1\} \prec\cdots \prec \{\sigma_1,\dots ,
\sigma_{2n}\} $.
If the numbers $1,\cdots ,2n$ are placed in a table, as
in Figure \ref{wreathbox}, then each element of $S_2\wr S_n$ 
may be viewed as the 
composition of some $\pi_1\in S_{2n}$ permuting the elements of each
row with a permutation $\pi_2$ permuting the rows.   
Thus, $\pi_1 = (12)^{e_1}
(34)^{e_2}\cdots (2n-1,2n)^{e_n}$ for some choice of $e_i\in \{0,1\}$ for
each $1\le i\le n$, and there is some $\pi\in S_n$, such that 
$\pi_2 (2i) = 2\pi (i)$ and $\pi_2(2i-1) = 2\pi(i)-1$ for $1\le i\le n$.

The action of $S_2\wr S_n$ on $\{ 1,\dots ,2n\} $ induces an action on the
saturated chains in $B_{2n}$.
Let us denote orbit representatives of this action by the
permutations in $S_{2n}$, written in one-line
notation, which label the chosen saturated chains in $B_{2n}$.
We choose the saturated chain labelled by the lexicographically 
smallest possible permutation in an orbit as the orbit 
representative.  The permutations in $S_{2n}$ occurring as labels are the
ones with the property that $2i-1$ comes before both $2i$ and $2i+1$ for
$1\le i<n-1$ 
and that $2n-1$ comes before $2n$.  
The labelling for orbits is a chain-labelling using the labels assigned to 
the orbit representative.  

\begin{examp}
\rm
The orbit representatives for $\Delta(B_6)/S_2\wr S_3 $, listed in
lexicographic order, are 
$123456$,
$1235\bullet 46$,
$123\circ 56\bullet 4$,
$13\bullet 2456$,
$13\bullet 25\bullet 46$,
$13\bullet 256\bullet 4$,
$1\circ 34\bullet 256$,
$1\circ 345\bullet 26$,
$1\circ 3456\bullet 2$,
$135\bullet 246$,
$13\circ 5\bullet 26\bullet 4$,
$1\circ 35\bullet 4\bullet 26$,
$1\circ 35\bullet 46\bullet 2$,
$13\circ 56\bullet 24$, and 
$1\circ 3\circ 56\bullet 4\bullet 2$.
Hollow dots denote swap ascents while filled-in
dots indicate descent locations.  For instance,
the swap ascent in $1\circ 3456\bullet 2$ comes from a
codimension one face skipping rank 1 in the intersection of
$134562$  with $132564$, resulting from
the fact that $312564$ is in the same orbit as $134562$.
\end{examp}

In the shelling argument for $\dwreath $, 
we refer to the $i$th row as being {\bf
empty} at the element $u$ in a saturated chain 
$\hat{0}\prec \{ \sigma_1 \} \prec \cdots \prec \{\sigma_1\cdots \sigma_k\} =u$ if 
$2i-1,2i\not\in \{ \sigma_1 ,\cdots ,\sigma_k \}$.  Similarly
we call row $i$ {\bf half-full} at $u$ if 
$2i-1\in \{ \sigma_1 ,\cdots ,\sigma_k \}$ but $2i\not\in
\{ \sigma_1 ,\cdots ,\sigma_k \}$, and we refer to row $i$ as {\bf full}
at $u$ if $2i-1,2i\in\{ \sigma_1 ,\cdots ,\sigma_k \}$.
 
All the saturated chain orbits belonging to the same rooted 
edge-interval $\hat{0}\prec u_1\prec \cdots
\prec u_k = u < v$ must agree in the following three ways:
\begin{enumerate}
\item
The same collection of half-full rows of $u$ must be full in $v$ 
\item
The same number of rows must switch from empty in $u$ to full in $v$.
\item
The same number of rows must switch from empty in $u$ to half-full in $v$. 
\end{enumerate}
The first of these three conditions depends on our use of 
edge-intervals which are rooted, since any two half-full rows in 
$u$ are equivalent, but they are distinguishable within the context of
a saturated chain orbit from $\hat{0} $ to $u$ .  In order to verify the
increasing-chain condition below, we will show that there cannot be two
(topologically) increasing chains that agree in all three ways.  

\begin{thm}
The labelling of saturated chains with minimal orbit representatives 
gives a lexicographic shelling.  More specifically, it is a $CC$-shelling.
\end{thm}

\proof
The crossing condition follows from Proposition \ref{quot}.
Next we classify topological ascents and descents in order to verify the
increasing chain condition.

We claim that every descent is an honest descent.  
Replacing a descent $\sigma_i \sigma_{i+1}$ 
by the ascent $\sigma_{i+1} \sigma_i$ 
yields a lexicographically smaller 
permutation, thus a member of a different orbit (since the orbit
representative was already the lexicographically smallest member of
its orbit).  In fact, the permutation obtained by this swap is the 
orbit representative of a new orbit 
since no new ascents have been introduced by swapping $\sigma_i $ and
$\sigma_{i+1} $, which means that
the requirement that $2j-1$ come before $2j$ and $2j+1$ for all $j$ is
still satisfied.

Whenever an orbit representative $F_j$ has
$2i-1$ immediately preceding $2i+1$
and then later has $2i+2$ before
$2i$, there is a 
swap ascent at the node between $2i-1$ and 
$2i+1$.  This is because the facet $F_i$ with $2i+2$ and 
$2i$ swapped is lexicographically earlier, and the facets $F_i, F_j$
share a 
codimension one face which omits the node between $2i-1$ and $2i+1$.
Notice that the characterization of this node as a swap ascent 
depends on later ranks in the chain.  This is not a 
problem, since the labels themselves, and thus the lexicographic order
on facets, only depends on the root of the chain.
It is not hard to check that all 
other ascents are honest ascents, and we have already shown that there
are no swap descents.

Now we must verify that each rooted interval has a unique (topologically) 
increasing chain.  If two orbit representatives 
both had only topological ascents on a rooted interval from $u$ to 
$v$, each would have labels in increasing order and be  
free of swap ascents of the type described above.  To avoid descents, 
each increasing chain must begin by 
completing all the requisite half-filled rows of $u$ in increasing
order, and then begin new rows, proceeding sequentially.  
The lexicographically smallest chain in the rooted
edge-interval proceeds through all the empty rows of $u$ to be filled in $v$
before turning to empty rows of $u$ to be half-filled in $v$.
Any lexicographically later orbit on the interval which is free of descents 
would also first complete the set of half-full rows of $u$ to be filled in 
$v$ and then would proceed to new rows.  To
differ from the lexicographically smallest orbit it would at some point 
necessarily insert one element into a new row immediately
before inserting two elements to the next new row.  
This would yield a swap ascent, as discussed above.  For example, 
the orbit labelled $134$ comes after the increasing 
chain labelled $123$ in the interval $(\emptyset ,\{1,2,3 \} ) $
in $B_4/S_2\wr S_2$, but $134$ has a swap ascent at $13$ since $4$ comes
before $2$ in any orbit including $134$.  Thus, the
lexicographically smallest orbit on the interval is the only topologically 
increasing orbit, just as desired.

Finally, we must confirm the multiple-face-overlap condition.
Suppose two faces $\sigma $ and $\tau $ are maximal in some 
$F_i\cap F_j$ and that $\sigma $ is maximal in $F_j \cap 
(\cup_{i'<j}F_{i'} )$.  
We may assume that $\tau $ has support $1,\dots ,r,s,\dots ,2n-1$,
because we will show (in the last paragraph of our proof) that for
$\Delta (B_{2n})/S_2\wr S_n$ every maximal face in $F_j \cap 
(\cup_{i'<j} F_{i'})$ has support omitting a single 
interval.  Assume also that $\sigma $ has support including 
$1,\dots ,r'$ for some $r'\ge s$ such that
no face in $F_i\cap F_j$ has support $1,\dots ,r'+1$.  
In particular, this means that $F_j$ 
agrees with $F_i$ on the interval skipped by $\tau $.
Let $u$ be the element of rank $r'$ in $F_j$, let $v$ be
the element of rank $r'+1$ in $F_j$ and let $v'$ be the 
element of rank $r'+1$ in $F_i$.  Notice that $u\prec v$
is equivalent to $u\prec v'$ in $\tau $, but not in
$\sigma $.  This is only possible if two rows are equivalent
in $\tau $ but not in $\sigma $, i.e. ~if two different rows
switch from empty to half-full on the interval skipped by
$\tau $.  Furthermore, the earlier of these rows (call it $\rho' $) must be
completed in the step $u\prec v'$ while the later row (denoted by $\rho $)
is completed in the step $u\prec v$.  We may assume that the 
interval skipped by $\tau $ has no topological descents,
for otherwise $\tau $ would be contained in a codimension one face of
$F_j \cap (\cup_{i'<j} F_{i'})$, and we would be done.  

Thus, all the letters inserted into rows to be half-filled
in the interval of $F_j$ skipped by $\tau $ are inserted by 
progressing through these rows sequentially.
If $\rho'$ is half-filled 
immediately before $\rho $, then there is a codimension one face
within $F_j$ that contains $\tau $ and belongs to $F_i\cap F_j$.
This face is obtained by half-filling $\rho' $ and $\rho $ in a 
single step.  Otherwise, let $\rho'=\rho_1 ,\rho_2,\dots ,\rho_k=\rho $ be 
the sequence of consecutive rows that are half filled between $\rho'$ and 
$\rho $ in $F_j$.
There must be some $\rho_{m-1},\rho_m$ for $1<m\le k$ such that 
$\rho_m $ is complete before $\rho_{m-1}$ in $F_j$, since $\rho $ is 
completed before $\rho'$ in $F_j$.  There also must be an $F_{i'}$ which 
agrees with $F_j$ except that it reverses the order in which the 
second elements of rows $\rho_{m-1}$ and $\rho_m$ are inserted.  This 
$F_{i'} $ comes before $F_j$ lexicographically.  There is a codimension
one face of $F_j$ that contains $\tau $ and is contained in 
$F_j \cap F_{i'}$.  It is obtained by 
half-filling $\rho_{m-1}$ and $\rho_m $ in a 
single step.  Thus, $\tau $ is contained in a codimension one face
belonging to $F_j \cap (\cup_{i'<j} F_{i'}) $, as desired.

There are no facets $F_j$ such that $F_j \cap (\cup_{i<j} F_i)$ has a maximal
face whose complement has disjoint support because this would mean that
skipping disjoint intervals allows two covering relations to be identified
that could not be identified by skipping a single interval.  Such 
identification could only come from making two rows interchangeable, but
because the rows have length two, it suffices to skip the single interval
beginning where the first row is half-filled and ending where the second
row is half-filled.  Thus, we
have checked all the necessary conditions for a lexicographic shelling.
\EOP

\medskip
The story is quite a bit different for 
$\Delta(B_{kn})/S_k \wr S_n$ when $k$ is greater than $2$.  
The increasing chain
condition fails for the lexicographic order on 
orbit representatives chosen to be lexicographically as small as
possible.  Consider the first four
facets $F_1 = 123456, F_2=124356,F_3=124536,F_4=124563$ in $\Delta (B_6)/
S_3\wr S_2$.  The intersection
$F_4 \cap (F_1 \cup F_2 \cup F_3) $ has two maximal faces $\sigma ,\tau $
where $\sigma $ is colored $1,2,3,4$ and $\tau $ is colored $4,5$, so 
$\tau $ has codimension greater than 1.  Here, $\sigma\subseteq F_4\cap
F_3$ and $\tau \subseteq F_4\cap F_3$.

Indeed, $\Delta(B_6)/S_3\wr 
S_2$ cannot possibly be shellable.  Molien's Theorem (cf. [St1]) together
with code given to us by Vic Reiner allowed us to determine the 
Hilbert series 
for the ring of invariants $k[\Delta(B_6)]^{S_3\wr S_2}$ with generators 
graded by poset rank; we obtained
$(1+q^2+q^3+2q^4+q^5+2q^6+q^7+q^8)(1/((1-q)(1-q^2)(1-q^3)(1-q^4)(1-q^5)(1-q^6)))$.
However, if $\Delta(B_6)/S_3\wr S_2$
had a shelling, then Theorem 6.2 of [GS] (stated in
Section \ref{invtsec}) would also yield this Hilbert series as follows.
The product $1/((1-q)(1-q^2)(1-q^3)(1-q^4)(1-q^5))$ accounts for polynomials
in $\theta_1,\theta_2,\theta_3,\theta_4$ and $\theta_5$, as discussed in
Section \ref{invtsec} while the numerator $1+q^2+q^3+2q^4+q^5+2q^6+q^7+q^8$
accounts for elements of the basic set given by a theorem of [GS],  
described in Section \ref{invtsec}.  That is, each shelling step $F_i$ with
minimal new face $G_i$ of support $S_i$ contributes $q^{d_i}$ to the
numerator where $d_i = \sum_{x\in S_i} x$.

In the case of $\Delta(B_6)/S_3\wr S_2$, there is no 
ordering on facets that would give rise to the necessary
collection of exponents $0,2,3,4,4,5,6,6,7,8$, so there cannot
be a shelling.  To see this, let us represent facets by
lexicographically smallest orbit representatives.  Note
that whichever facet comes last among $142536,142563,145236$ and $145263$
will contribute a minimal new face whose support includes ranks
1,3 and 5.  This implies that the last step among these four  
would contribute an exponent of 
at least 9 to the numerator of the 
Hilbert series, a contradiction since $q^8$ is the largest
power of $q$ present.

\section{A lexicographic order on facets of $\Delta (\Pi_n )/S_n $}
\label{partsec}

First we briefly describe which chains in $\Pi_n$ are in the same orbit
and then give a way of representing chain orbits by trees.  After this,
we give a chain-labelling on $\dps $ in terms of this tree-representation
for the facets of $\dps $.  We remark that a similar tree-representation 
appears in [Ko3].

\subsection{Orbits of partition lattice chains}
The $S_n$-orbit of a partition in $\Pi_n $ is 
determined by the size of its blocks, i.e. by an integer partition.  Thus, 
vertices in $\dps $ are given by unordered partitions of the 
integer $n$.  However, 
two edges $u < v$ and $u' < v'$ in $\Delta(\Pi_n) $ may belong to distinct
$S_n$-orbits even if $u'\in Orb(u)$
and $v'\in Orb(v)$, because this does not guarantee the existence
of a single permutation
$\pi\in S_n$ such that $\pi(u) = u'$ and
$\pi(v) = v'$.  If no such $\pi $ exists, then
$u<v$ and $u'<v'$ give rise to distinct edges in
$\dps $ which have the same vertices.
Thus, faces in $\dps $ are 
sequences of successively refined partitions of the integer $n$ 
enriched with some additional information.

\begin{examp}
\rm
Consider a chain $22\prec 11|11\prec 3|8|11$ with numbers denoting block
sizes.  The orbit does not depend on 
which block of size 11 is split into blocks of size $3,8$.
However, there are two orbits of the form 
$22\prec 11|11 \prec 3|8|11 \prec 3|8|3|8 \prec 3|4|4|3|8$, because the two 
blocks of size 8 are not interchangeable since they were created at 
different ranks.
On the other hand, the two blocks of size 8 are indistinguishable within 
the chain which skips immediately from $11|11$ to $3|8|3|8$.
\end{examp}

The orbit of a chain keeps track of
what type of block is split into what types of pieces at each rank in a
chain.  To be precise, the {\bf type} of a block in $u_k$ will 
be its $S_n$-equivalence class, relative to a chain 
$\hat{0}<u_1<\cdots < u_k < \hat{1}$.  Two blocks $b_1,b_2$ in $u_k$ 
are said to be {\bf $S_n$-equivalent} if (1) they have the same content, 
(2) they must be created in the same refinement
step $u_i < u_{i+1}$, and (3) they must be children of 
blocks $B_1$ and $B_2$, respectively, which are 
themselves equivalent to each other.  
This includes the possibility that $B_1=B_2$, so $b_1,b_2$ come from a 
single parent.  The point is that $S_n$ may swap the elements of $b_1$ with
the elements of $b_2$ in a way that preserves the chain.  
Denote the orbit of a chain $C$ by $\pi(C) $.

\subsection{Two encodings for faces in $\dps $}
First we encode the orbit of a chain as a 
tree whose nodes are the partition blocks that occur in the chain.
The root is the single block appearing as $\hat{0}$ in any poset chain.
The children of a tree node $B$ are the blocks obtained from $B$ when it is
refined in the chain.  Each tree node is labelled
by its block content, and each parent is also labelled with  
the rank at which it is refined.  The equivalence class of 
a particular block relative to a chain orbit comes from the 
tree built by that chain orbit.  
Two blocks are equivalent if there is a graph automorphism that swaps
the blocks and carries each tree block to one with identical labels.
The orbit of a saturated chain 
gives rise to a binary tree.

%\begin{figure}[h]
%\begin{picture}(175,150)(0,0)
%\psfig{figure=hist.eps,height=5cm,width=8cm}
%\end{picture}
%\caption{Equivalent blocks in terms of trees}
%\label{history}
%\end{figure}
%In Figure \ref{history}, the nodes labelled
%A are equivalent (assuming corresponding node labels of size and 
%rank agree),
%but the node labelled B is not equivalent to either node labelled A.

The labelling given next will depend on a choice of planar embedding
for these trees.  We will define this embedding by specifying for each
parent an ordering on its children.  
Let us make a convention for choosing a planar embedding so as to 
assign a label $\lambda (\pi (\hat{0}\prec\cdots \prec u\prec v)) $ to each covering 
relation orbit $\pi (u\prec v)$ based on the entire saturated chain orbit 
$\pi (\hat{0}\prec \cdots \prec u)$ to which $\pi (u\prec v) $ belongs.
In the process of choosing the label, we shall also choose a planar 
embedding for the tree given by $\pi (\hat{0}\prec \cdots \prec u\prec v)$.  
In particular, this embedding will give us an ordering on the blocks for
each partition in the chain, which depends only on the root of that chain.
The chain-labelling, provided next, uses this ordering on the
blocks of a partition to be refined in a covering relation.

Sometimes, we will use a more compact encoding for a chain.  Namely, we 
list $n$ balls in a row with $n-1$ bars separating them, and place numbers
between $1$ and $n-1$ below the bars.  The balls represent the $n$ numbers
being partitioned, since these are freely interchangeable.  The numbers 
below the bars record ranks at which bars are inserted while progressively
refining a partition.  This representation for a chain often is not unique,
but it completely determines the chain orbit.  We will order the blocks
in the fashion described at the beginning of Section ~\ref{nearicc}, so this
will determine our choice of representation.  See Figure ~\ref{twoblock} for 
an example.

\subsection{A chain-labelling which nearly satisfies the increasing chain
condition}\label{nearicc}

Let us assign labels to orbits
$\pi (\hat{0}\prec u_1\prec \cdots \prec u_k = u \prec v)$ of rooted
covering relations.
The distinct $\pi (u\prec v)$ with a fixed root
$\pi (\hat{0}\prec u_1 \prec \cdots \prec u_k = u )$ are specified by 
which type of block from $u$ 
is split (recalling that type means equivalence 
class) together with the content of its children.
Assume by induction that a saturated chain orbit $\pi (\hat{0}\prec 
u_1\prec\cdots \prec u_k=u )$ imposes an order on the blocks of $u$.
We obtain from this an ordering on the blocks of $v$ for each $u\prec v$
as follows:
\begin{enumerate}
\item 
Split the leftmost block $B$ in the ordered partition of $u$ which belongs to the 
equivalence class to be split by the covering relation $u\prec v$.
\item
Place the lexicographically smaller of the two blocks obtained from $B$ 
to the left of an inserted bar separating the two blocks resulting from $B$.
\end{enumerate}
Thus, we get an ordering on the blocks of $v$ from the ordering for $u$ by 
replacing $B$ by the two
blocks derived from it ordered lexicographically and otherwise preserving
block order.  One may then use the
position where the new bar is inserted to give a chain-labelling.

\begin{rk}
{\rm
This labelling by bar position is motivated by a feature of $\Pi_n $
(which is related to the splitting basis for the partition lattice
given in [Wa]).
Each permutation $\pi\in S_n$ gives rise to a boolean sublattice
$B_{n-1} $ of those partitions obtained by 
listing $1,\dots ,n$ in the order given
by $\pi $ (written in one-line notation).  Each partition consistent with 
$\pi $ is specified by choosing a subset of the $n-1$ possible bars to
insert splitting the numbers into blocks. 
}
\end{rk}
 
The labelling by bar position does not always
satisfy the increasing
chain condition, as indicated by Example \ref{noshell}.  
To transform this labelling into one which will 
satisfy the increasing chain condition,
we will introduce a block-sorting step next.

\begin{examp}\label{noshell}
\rm
Consider $\Delta (\Pi_{22} )/S_{22} $.  Take the rooted
interval $\hat{0}\prec u < v$ with 
$u = 11|11 $
and $v=1|1|9|2|2|7$ where one block of size 11 splits into 
$1|1|9$ and the other into $2|2|7$.  
This rooted interval includes the
product of chains in Figure \ref{elevenint} with chain-labels as shown.
The left half gives labels for saturated chains beginning with $11|11 \prec 1|10|11$
and the right half for those beginning with $11|11\prec 2|9|11$.  
The lexicographically smallest chain 
appears farthest to the left and is labelled $1,2,13,15$.  The rightmost
chain in Figure \ref{elevenint} 
consists entirely of honest ascents, even in light of 
all the saturated chains in the interval rather than only those depicted
in Figure \ref{elevenint}.  
Thus, this chain-labelling violates the 
increasing chain condition on the interval from $11|11$ to $1|1|9|2|2|7$.  
\begin{figure}[h]
\begin{picture}(200,150)(45,0)
\psfig{figure=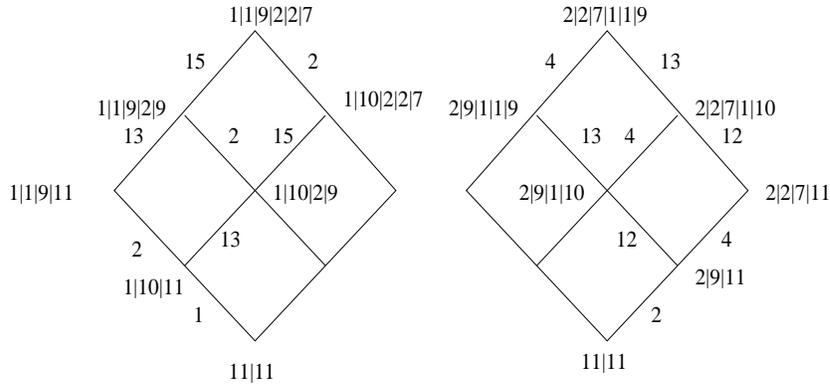,height=5cm,width=11cm}
\end{picture}
\caption{First approximation at a chain-labelling}
\label{elevenint}
\end{figure}
Figure \ref{elevenint} shows
the chains in the interval in which bars are inserted left to right.
Covering relations are labelled by bar positions.  However, the labelling
does not give a shelling
because the intersection of the rightmost chain with those coming 
earlier is not pure of codimension one.  
\end{examp}

\subsection{A modified chain-labelling which satisfies the increasing chain 
condition}

\begin{figure}[h]
\begin{picture}(100,120)(0,0)
\psfig{figure=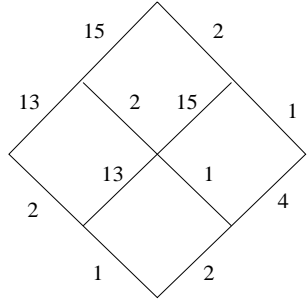,height=4cm,width=4cm}
\end{picture}
\caption{Revised chain-labelling which incorporates block sort}
\label{goodelvs}
\end{figure}

We will add a block sorting step just before bar
position is recorded.  This will yield the labelling in
Figure \ref{goodelvs} for the product of chains from Figure
\ref{elevenint}.  
When we split a block of
size 11 into blocks of size 2,9, and then break the other block
of size 11 into blocks 1,10, we sort the blocks of size
1,10 to the left of the blocks of size 2,9 before assigning the second
label.  

Our chain-labelling assumes (by induction) that we have already
assigned a label to the orbit $\pi (\hat{0} \prec
u_1\prec\cdots\prec u_k=u)$ and in the process have
ordered the blocks of $u$.  It then specifies a label for $u\prec v$ as
well as an ordering on the blocks of $v$.

\begin{enumerate}
\item
Block refinement: Split the leftmost block $B$ belonging to the  
equivalence class to be split in $u\prec v$ 
into two blocks, the smaller of which is placed on the left.  Thus,
a block $B$ is replaced by two blocks derived from it with the smaller
on the left, giving a more refined ordered partition.
Note that block equivalence in orbits of saturated chains 
only comes from a single block
splitting into two identical pieces.
\item
Block sort: If $B$ is equivalent to another block $B'$ in $u_i$ for 
some $i<k$, let $P(B)$ be the common parent of $B$ and $B'$.
Compare the subtree of descendents of $B$ to the subtree
of descendents of $B'$ to decide whether $B$ or $B'$ should be the left 
child of $P(B)$ in the block order for $v$.  The left child is the
block with the lexicographically smaller word comprised of the positions
at which bars occur within that block.
Ties are broken using the rank at 
which $B$ and $B'$ were first split (the block with earlier rank is 
sorted farther to the left).  Next, apply this sorting procedure 
to the successive ancestors of $B$.
\item
Chain-labelling: The label assigned to 
$\pi (\hat{0}\prec\cdots \prec u\prec v)$ is
a 3-tuple with the following components, listed in order of precedence:
\begin{enumerate}
\item
The post-sort position of the newly inserted bar.
\item  
The word consisting of all the post-sort bar positions in $v$.
\item
The (ordered) list of ranks at which the successive 
ancestors of $B$ were themselves split, with $P(B)$ given lowest
precedence.  
\end{enumerate}
\end{enumerate}

By Proposition \ref{quot}, $\Delta (\Pi_n)/S_n$ satisfies the crossing
condition, so let us turn our attention
to the increasing chain condition and multiple-face-overlap condition.
The proofs of Theorem ~\ref{icc} and Theorem ~\ref{mult-face}
adopt several ideas from [HK].

\begin{thm}\label{icc}
The above labelling satifies the increasing chain condition.
\end{thm}
\proof
Denote the above labelling 
by $\lambda $.  Notice that any topologically increasing chain on a rooted
interval $\hat{0}\prec\cdots\prec u < v$ must insert bars from left to right 
(post-sort) into $u$ in such a way that each block of $u$ is split into 
pieces which are nondecreasing in size from left to right.  There is at least
one such increasing chain, obtained by greedily inserting bars from left to 
right so that the children of each block are nondecreasing in size from left 
to right.  Whenever two equivalent blocks are to be refined, refine them in
such a way that the second one to be refined does not get sorted to the left
of the first one refined.  That is, refine the first of these two blocks in
the way that gives a smaller left child as soon as the two refinements differ.
What remains to show is that 
there is a unique topologically increasing chain on each interval and that 
this chain is lexicographically earliest on
the interval.  To this end, we will show that any saturated chain that
is not lexicographically first on an interval has a topological descent.  Our 
argument will repeatedly use the fact that sorting never moves newly 
inserted bars to the right.

Suppose that the orbit $\pi(C)$ of a saturated chain 
$\hat{0}\prec u_1\prec\cdots \prec u_k\prec\cdots\prec w$
is not lexicographically smallest on a rooted edge-interval 
$\hat{0}\prec u_1\prec\cdots\prec u_k < w$. 
Consider $u\in \pi(C)$ of lowest rank such that the rooted 
edge $\pi (\hat{0}\prec\cdots\prec u\prec v)$ in $\pi(C)$
has larger chain-label 
than $\pi (\hat{0}\prec\cdots\prec u\prec v')$ for some $v'\ne v$ 
in a different saturated chain orbit also
belonging to the rooted edge-interval
$\pi (\hat{0}\prec\cdots\prec u_k < w)$.  
Choose $v'$ with the minimal possible label 
among all such choices.  Let $B$ be the block from $u$ that is split in the 
step $u\prec v$ and let $B'$ be the block (also from $u$) that is
split in $u\prec v'$.  We consider three cases, depending on whether
the labels $\lambda (u\prec v),\lambda(u\prec v')$ 
first differ in the first, second or third component
of the label.  We will sometimes abuse notation by not listing the root,
even though the labellings do depend on them.

Case I: Suppose the labels differ in the (post-sort) bar position. 
We further subdivide this case, depending 
on whether the blocks $B,B'$ being split by $u\prec v$ and $u\prec v'$, 
respectively, are equal or not.  Case Ia: If $B=B'$, then
we will show that $\pi(C)$ must either have an honest
descent or a swap ascent within the interval.  Let $b'$ be the 
smaller of the two blocks resulting from splitting $B'$ in 
$u\prec v'$; note that $b'$ must eventually be derived from
$B$ or an equivalent block later in $\pi (C)$, within the interval (since 
$u\prec v'$ belongs to the same interval).
Before this happens, there will be a descent or swap ascent since this 
later step will have a lexicographically smaller label than any step before
it in the interval.  This is because the block $B$ (or an equivalent block
from which $b'$ is obtained) is at this point sorted at least as far to the
left as it would be when the bar is inserted in $u\prec v'$.
Thus, the bar splitting off $b'$
must be to the left of the bar insertion immediately before it
(and thus there must be a descent) 
unless there is another step splitting $B$ into two larger pieces
earlier in $\pi (C)$ such that $b'$ comes from the right component
in this split block.  In this case,
there must either be a swap ascent in $B$ immediately before $b'$ is split off
on its own, or there must be a descent if
there are steps splitting $B$, then splitting another block, then later
splitting a block derived from $B$ to create $b'$.

Case Ib: If $B \ne B'$, then eventually we will split $B'$.
At this point, either we break off the smallest piece of $B'$, just as 
in $v'$, yielding the smallest label so far in $\pi(C)$ restricted to the 
interval and thus 
a descent, or else the preceding $B = B'$ argument may be applied to the 
remainder of $\pi(C)$ within the interval to obtain a topological descent at
a higher rank.  

Case II: Now suppose $u\prec v'$ has the same post-sort
bar position as $u\prec v$, but that their labels
differ in the words made up of post-sort bar positions in $v$ and $v'$.  
Then $B\ne B'$ when we let $B, B'$ be the blocks which are split in
$u\prec v$, $u\prec v'$, respectively.
Eventually, a bar must be inserted in $\pi(C)$ into a 
block which is equivalent to $B'$ to create a smallest possible
piece $b'$, as above.  Note that 
equivalent blocks in a saturated chain orbit must be
adjacent regardless of how they are split and sorted, by virtue of the 
underlying binary tree of blocks.  Therefore, when $b'$ is derived from $B'$
later in $\pi(C)$, sorting will still move the
new bar to the position it would have achieved in $v'$, 
unless $B'$ is split into two 
larger pieces and then $b'$ is derived later from the larger of these.  
When we are not
in this special case, then the bar position
word when $b'$ is created will be smaller than that for 
$u\prec v$, and the bar to the 
right of $b'$ is placed at least as far to the left as it would be in 
$u\prec v'$.  Thus, 
$\pi (C) $ must have a descent on the interval.  In the special case 
where $b'$ is 
derived from a descendent of $B'$ which is not the leftmost descendent, 
there still 
will be a swap ascent or a descent.  The former must occur if one of the 
steps splitting 
$B'$ occurs immediately before $b'$ is split off.  The latter must occur 
if $B'$ is
split, then some $B $ not descending from $B'$ is split and then later $B'$ is 
further refined.

Case III: Now suppose that the first two coordinates agree, 
but that the earliest distinct ancestor of $B'$ is split at an 
earlier rank than that of $B$.  This means 
$B$ and $B'$ are not equivalent.  Hence, blocks of
both types must be split in any saturated chain orbit in the interval, 
yielding a 
descent or swap ascent in $\pi(C)$ eventually, since the block $B'$ is 
sorted to the
left of $B$ when it is split.  If $B'$ is split immediately after 
$u\prec v$, then there is
a lexicographically smaller chain with these two steps reversed, i.e. a 
topological descent.
Otherwise, any intermediate steps will have larger labels, and the 
chain must have a descent (or swap ascent) by similar reasoning to above.
Thus, we have verified the increasing chain condition in all cases.
\EOP

\begin{cor}
The quotient complex $\Delta(\Pi_n)/S_n$ has lexicographic shelling
steps at all facet insertions except those which violate the
multiple-face overlap condition.
\end{cor}

Next we describe precisely where the multiple-face-overlap condition 
fails.

\begin{thm}\label{mult-face}
The above labelling also satisfies the
multiple-face-overlap condition, except in the case of
certain facets with identical blocks
created in consecutive steps from a single parent.
\end{thm}

\proof
We will show that the multiple-face-overlap condition holds 
assuming that the two blocks to be identified have the same parent block.  
Suppose 
$F_j\cap F_k$ has two maximal faces $\sigma ,\tau $ such that $\sigma $ is 
maximal in $F_k \cap (\cup_{i<k} F_i)$ and $\tau $ has support including
$1,\dots ,r'$ where no face in $F_j\cap F_k$ has support $1,\dots ,r'+1$. 
We need only check the condition if 
$\sigma $ has support $1,\dots ,r,s,\dots ,n$ for some $s<r'$.
We may also assume $F_k$ consists
entirely of topological ascents on the interval from rank $r$ to $s$.
Otherwise, $\sigma $ has codimension one, and we would be done.  
Thus, the post-sort positions of the bars inserted within $F_k$ from 
rank $r$ to $s$ are increasing left to right.  Let $u$ be the element of rank
$r$ in $F_k$.  Then the bars inserted into any particular block
of $u$ between rank $r$ and $s$
are also arranged so that blocks are nondecreasing in
size from left to right.  

Note that the splitting step between ranks $r'$ and $r'+1$ in $F_j$
must not be equivalent to the splitting step in $F_k$.
However, when we restrict to $\sigma $, these steps
become equivalent relative to the chain $\sigma $.  
Thus, the blocks being split must have the same
content.  Let us call these blocks $B_j$ and $B_k$, respectively.  If 
$B_j$ and $B_k$ are created in consecutive steps in the interval of 
$F_k$ from rank $r$ to $s$, then we get a codimension one face belonging to 
$F_j\cap F_k$ 
which contains $\sigma $.  This face is obtained by skipping the step
immediately after $B_j$ is created and immediately before $B_k$, so that 
these are
created in a single step from the same block.
If $B_j$ and $B_k$ are not created consecutively in $F_k$ but do have the
same parent, then
there is some block $B_{k-1}$ created immediately before $B_k$ in 
$F_k$ which is
also equivalent to $B_k$ in $\sigma $.  Consider the facet $F_{j'}$ with
$j'<k$ which agrees with $F_k$ 
except that the roles of $B_k$ and $B_{k-1}$ are reversed from rank $r'$
onward.  $F_{j'}$ comes before $F_k$ lexicographically because the covering
relations where they first differ share the same first two coordinates of 
their label.  Because $B_{k-1}$ was created earlier, the third 
component of the label is smaller for $F_{j'}$.  
Since $B_{k-1}$ and $B_k$ are created in consecutive steps in
$F_k$ in the interval skipped by $\sigma $, the codimension one face 
creating $B_{k-1}$ and $B_k$ in a single step in $F_k$ 
will be in the intersection $F_k \cap (\cup_{k'<k} F_{k'} ) $ 
and will contain $\sigma $, as desired.  
\EOP

Section ~\ref{partsec} will discuss what to do when
$B_j$ and $B_k$ come from distinct parent blocks.
It will use a partitioning for real projective space and some 
generalizatinos to incorporate non-shelling steps 
into the partitioning.  First, we give an example of a face in 
$\Delta (\Pi_8)/S_8)$ whose link is $\reals P_2$, implying for $n\ge 8$
that $\Delta (\Pi_n)/S_n$ is not Cohen-Macaulay over $\integ /2\integ$,
and hence is not shellable.  The above
analysis helped us find this example.

\begin{examp}
\rm{Consider the chain 
$$12345678 < 1234|5678 < 1|234|5|678 < 1|2|34|5|6|78 < 1|2|3|4|5|6|7|8$$
in $\Pi_8$.  Its $S_8$-orbit is a face in $\Delta (\Pi_8)/S_8$ whose link
is the real projective plane.  To see this, notice that the link has three
vertices: one at rank 1, one at rank 3 and one at rank 5, and it has 6
edges and 4 2-simplices, all arranged as $\reals P_2$.  This is 
essentially the same construction used in Proposition 3.1 of [He]
to study $\Delta (B_{lm})/S_l \wr S_m$, so more detail may be found there.}
\end{examp}

Notice that our labelling never gives any chains 
consisting entirely of topological descents.
This is just as one would expect, given that $b_{1,\dots ,n-2} (n)$
vanishes [Ha], [St2] and that $\dps $ is collapsible [Ko1].  

\section{Partitioning $\Delta (\Pi_n )/S_n$}

\begin{figure}[h]
\begin{picture}(130,40)(45,12)
\psfig{figure=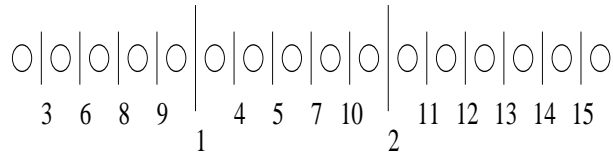,height=2cm,width=8cm}
\end{picture}
\caption{Facet involving two equivalent blocks}
\label{twoblock}
\end{figure}

Theorem ~\ref{icc} showed that the only way a non-shelling 
step $F_j$ may arise is when the multiple-face-overlap condition is
violated.  More specifically, Theorem ~\ref{mult-face} showed that
the facet $F_j$ must at some rank refine a block which in the context
of a face $F\subseteq F_j$ is equivalent to a block to its left.
For example, consider the face of support $2,4$ within the 
facet shown in Figure ~\ref{twoblock}.

\begin{defn}
A {\bf partitioning} of a pure boolean cell complex $\Delta $
is an assignment to each facet $F_i$ of one of its faces $G_i$, 
so that the boolean upper intervals $[G_i,F_i]$ give
a partitioning of the faces
in the complex, i.e., so that $\Delta $ may be written as
a disjoint union $[G_1,F_1]\cup\cdots\cup [G_s,F_s]$.  
\end{defn}

A partitioning of a pure, balanced boolean cell
complex $\Delta $ gives a combinatorial
interpretation for each flag $h$-vector coordinate $h_S(\Delta )$ as the 
number of facets $F_i$ whose minimal face $G_i$ has support $S$.  
In this section, 
we give a partitioning for $\dps $ as follows.  We show how to partition
subcomplexes of links of 
faces coming from progressively more general classes of
non-shelling steps.  Then we use the topological
ascents and descents of the lexicographic
order of Section 4 to extend and merge these link subcomplex 
partitionings into a partitioning for $\dps $.  

En route to a partitioning for $\Delta(\Pi_n)/S_n$, 
Section ~\ref{linkpart} 
will give partitionings for subcomplexes of certain links.  For
example, consider the link of the face $\tau $ of 
nsupport $2,4,6,8,10,11,12,13,14,15$ contained in the 
facet depicted in Figure ~\ref{twoblock}.  This
facet will contribute to our partitioning an interval 
$[G_j,F_j]$ with $G_j$ of support $2,5,9$.

\subsection{Partitioning subcomplexes of links arising in non-shelling steps}\label{linkpart}
Consider the face $F$ 
comprised of the following elements.  First it has a partition into 
$m$ blocks of size $k+1$ along with a single block of size $l$ for 
some $l>k+1$.  Then, for each $1\le j\le k$, the chain includes a partition
in which the block of size $l$ is unrefined and each of the other $m$ blocks is
split into $j$ singletons and a single block of size $k+1-j$.
Finally, the chain includes a partition whose only nontrivial nontrivial block
is the block of size $l$, and a saturated chain upward from this rank which
sequentially splits off singletons from the lone nontrivial block.
Let $\Delta $ denote the subcomplex of the link of $F$  
consisting of those chains which initially insert from left to right the
bars separating the $m$ equivalent blocks.

\begin{figure}[h]
\begin{picture}(120,40)(50,12)
\psfig{figure=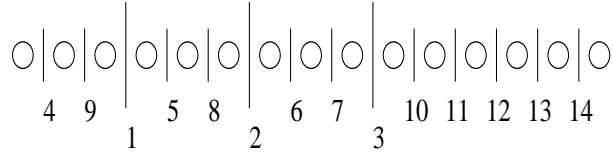,height=2cm,width=8cm}
\end{picture}
\caption{Facet with slot permutations $123,123,321$}
\label{threeblock}
\end{figure}

The facets in $\Delta $ correspond to $(k+1)$-tuples $(\sigma_0,\dots ,
\sigma_k)$ of permutations in
$S_m$ with the identity as the first permutation.  Requiring $\sigma_0$ to 
equal the identity reflects our requirement that the $m$ equivalent blocks
of size $k+1$ be created by inserting bars from left to right.
The permutation $\sigma_i\in S_m$ for $i>0$
specifies the order in which bars are inserted into the $i$-th slot in the 
$m$ blocks.  By $i$-th slot, we mean the bar position
separating $i$ objects to its left from $k+1-i$ objects to its right.
For example, the facet in Figure ~\ref{twoblock} has permutations
$12, 12, 21, 21, 12$ and the facet of Figure ~\ref{threeblock} is given 
by permutations $123$, $123$ 
and $321$.

\begin{defn}
We refer to the
$k$ positions into which bars may be inserted in the $m$ 
blocks as the {\bf slots}, and we call the
ranks at which one of the $m$ blocks is refined {\bf slot ranks}.  
Let us refer to the $m-1$ bars inserted left to right initially
to separate the $m$ 
equivalent blocks as the {\bf splitters}.  We say that a face includes a 
splitter if it includes the rank at which that splitter is inserted.
\end{defn}

\begin{figure}[h]
\begin{picture}(125,160)(110,7)
\psfig{figure=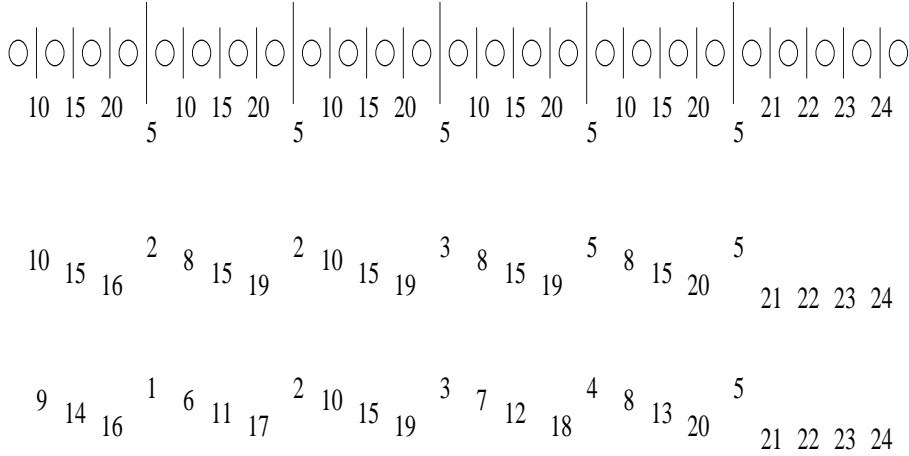,height=6cm,width=12cm}
\end{picture}
\caption{Face representation and extension}
\label{facerep}
\end{figure}

The facet given by $k+1$ copies of the identity permutation will have the
empty set as its partitioning minimal face.  For other facets, we are 
interested in
finding minimal rank sets needed to differentiate them from this
``identity facet''.
\begin{rk}
The effect of including a splitter in a face is to distinguish blocks to its
left from blocks to its right.  Including slot ranks from two different 
slots may record the fact that two blocks are split in opposite
order in the two slots.  For example, ranks 8,16 in Figure ~\ref{facerep}
show that blocks 1 and 2 are filled in opposite orders in the first and 
third slots.
\end{rk}
Figure ~\ref{threeblock} gives an example of a facet which does not have the
same face of support 1,8 as the identity facet.
Now let us specify how to assign minimal 
faces $G_i$ to the facets $F_i$, using the representation of facets as 
$(k+1)$-tuples of permutations $\sigma_0,\dots ,\sigma_k
\in S_m$.  
\begin{partcon}
{\rm
First we construct permutations $\pi_0,\dots ,
\pi_k$ from a facet $F_i$.  Let $\pi_{i+1}$ be the permutation
in two-line notation which has $\sigma_i$ (written in one-line notation) 
as the first line and $\sigma_{i+1}$ (again in one-line notation) as the 
second line.  The ``wrap-around'' permutation $\pi_0 $ is obtained by using
$\sigma_k$ as the first line and $\sigma_0$ as the second line.  The
minimal face $G_i$ associated to $F_i$ consists of the ranks $lm+j$ such that
$\pi_l^{-1}(j)> \pi_l^{-1}(j+1)$.  }
\end{partcon}

\begin{figure}[h]
\begin{picture}(100,120)(65,5)
\psfig{figure=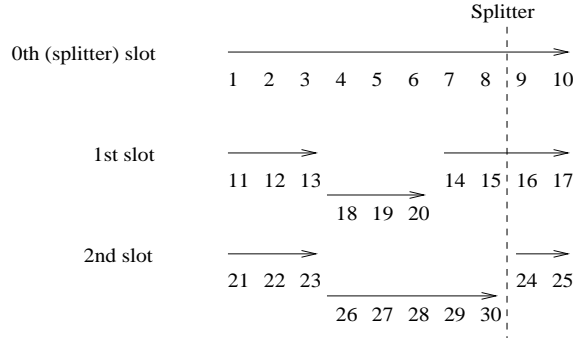,height=4.5cm,width=7.5cm}
\end{picture}
\caption{Facet with minimal face of support $8,17,28$}
\label{lategap}
\end{figure}

\begin{examp}
Letting $k=2, m=10$, Figure ~\ref{lategap} depicts the facet $F_i$ given
by permutations $\sigma_0=1,2,3,4,5,6,7,8,9,10$,
$ \sigma_1 = 1,2,3,7,8,9,10,4,5,6$ and $\sigma_2=1,2,3,9,10,4,5,6,7,8$.  
The minimal face assigned to $F_i$ has support $8,17,28$.
\end{examp}

To ensure that our assignment of minimal faces to facets
gives a partitioning, we must check (1) that every 
face belongs to some interval $[G_i,F_i]$ and (2)
that no face is included in two different intervals.  To verify (1), 
we describe how to extend any face $F$ to a facet $F_i$ whose minimal
face $G_i$ is contained in $F$.

\begin{prop}\label{cover}
Every face $F$ is contained in at least one interval $[G_i,F_i]$.
\end{prop}
\proof
First let us choose a representation for $F$, i.e. a specification of
which blocks (among equivalent choices) to split
to what extent at each slot rank in $F$.  
Figure ~\ref{facerep} depicts our choice for a  
face of support $2,3,8,16,19$.
By convention, we split blocks in a single slot as follows.  First, we
split equivalent blocks from left to right in the last slot.  At this
first step, two blocks are equivalent if there are no splitters separating
them.  
Now split each set of equivalent blocks in the penultimate
slot from left to right.  At this second step, two blocks must also 
have been split at the same rank in the last slot in order to be 
equivalent.  Proceed in this fashion through all the slots in
reverse order to obtain a representation of $F$.
Thus, $F$ is encoded as a choice of which actual blocks are split at 
each slot rank of $F$, along with the list of which splitter ranks 
belong to $F$.  

Now $F$ is extended to a facet $\overline{F}$ by refining the order in 
which bars are inserted as follows:
extend each interval between ranks in the first slot by splitting blocks
from left to right.  This choice determines the 
permutation $\sigma_1 $.  Similarly refine slot $s+1$, but instead of 
extending each interval from left to right, proceed in the order that is
increasing with respect to the permutation $\sigma_s$. 
This ensures that $\pi_{s}^{-1}$ will not have any descents on the intervals 
between consecutive slot ranks.  We may, however, have some
descents in $\pi_0^{-1}$ at ranks that are not splitter ranks.  

These ranks may be eliminated by modifying our face representation as 
follows.  Whenever 
$\pi_0^{-1}$ has a descent at a rank that is not a splitter rank,
this means there must be a slot rank at which two consecutive
blocks separated by this splitter cease to be equivalent.  We reverse the 
roles of these two blocks in our modified face representation.  The effect
is to change the descent in $\pi_0^{-1}$ to an ascent, and instead to have
a descent at the slot rank where the blocks cease to be equivalent.  Now
taking the increasing extension of this new face representation yields a 
facet $F_i$ such that the support of $G_i$ is contained in the support of
$F$, implying $F\in [G_i,F_i]$.  A very similar argument, with more detail 
included, is used to partition $\Delta(B_{lm})/S_l\wr S_m$ in [He].  
\EOP

\begin{examp}
Consider the link subcomplex $\Delta $ specified by the face listed at the
top in Figure ~\ref{facerep}.  Below this is our 
representation for a face $F$ of 
support $2,3,8,16,19$ in $\Delta $ and then its extension to a facet $F_i$. 
We list the ranks at which bars are inserted, using slight 
variation in height to distinguish different slots.
Observe from the descents in $\pi_0^{-1},\pi_1^{-1},\pi_2^{-1}$
that the minimal face $G_i$ has support $3,8,16,19$.
\end{examp}
Next we check that each face is included only once in the partitioning.  
\begin{prop}
There is no overlap among the intervals $[G_i,F_i]$.
\end{prop}
\proof
Each representation of a face $F$ has a unique extension to a facet 
that avoids descents between slot ranks in $F$.
Only one of these representations will also avoid descents from 
wrap-around at ranks that are not splitter ranks, making our choice
unique.  Thus, there is only one extension of $F$ to a facet $F_j$ such
that the support of $F$ contains the support of $G_j$, as needed.
\EOP

\begin{qn}
Is the incidence matrix for this matching of minimal faces with facets
containing them nonsingular?  If more generally this holds for a 
partitioning of $\dps $, then 
the partitioning will give a ring invariant basic set for 
$k[\Delta(\Pi_n)]^{S_n}$, by results of [GS].
\end{qn}

\subsection{Partitioning subcomplexes of larger links}\label{nextlink}
\begin{figure}[h]
\begin{picture}(180,105)(20,10)
\psfig{figure=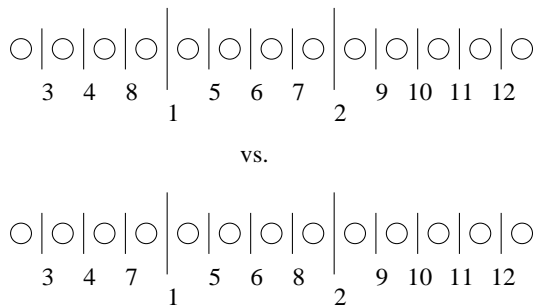,height=4cm,width=7cm}
\end{picture}
\caption{Two facets intersecting in 
faces of codimension one and three}
\label{j2}
\end{figure}

Next we generalize the partitioning of the previous section
to allow for examples such as the following:

\begin{examp}
Let $F_j, F_{i'}$ be the top and bottom facet
in Figure ~\ref{j2}, respectively.
$F_j \cap (\cup_{i<j} F_i )$ has a maximal face $G$ of support
$2,6,7,8,9,10,11,12$ since $G\in F_{i'}$ 
Consequently, $F_j$ would need to contribute
minimal faces of support $\{2,3,7\} ,\{2,4,7\} $ and $\{ 2,5,7\} $
in a lexicographic shelling.  
This is resolved in a partitioning by assigning the faces of
support $\{ 2,3,7\} $ and $\{ 2,5,7\} $ to lexicographically later 
facets.
\end{examp}

Consider a face consisting of the following chain 
$u_1 < u_2\prec u_3\prec\cdots\prec u_r $ of partitions.  Let $u_1$ consist
of $m$ equivalent blocks of size
$k+1$ and one block of size $l>k+1$.  Let $u_2$ have
these $m$ blocks refined to singletons and leave the block of size $l$
unrefined.  Finally, take a saturated chain 
$u_2\prec u_3\prec \cdots \prec u_r$ which at each stage splits off a 
singleton from its unique nontrivial block.
For example, take the face of support $2,8,9,10,11,12$ 
in either facet of 
Figure ~\ref{j2}.  Now we will consider the subcomplex of the link of such
a face in which we require that splitters be inserted strictly from left to 
right.  Denote this subcomplex by $\Delta $.

Each facet in $\Delta $
inserts splitters left to right, then refines the 
blocks $B_1,\dots ,B_m$ entirely.  
Our assignment of minimal 
faces $G_j$ to facets $F_j$ will generalize 
the partitioning of Section ~\ref{linkpart}, 
but it will use a single evolving permutation to play the roles of the 
permutations $\sigma_0,\dots ,\sigma_k $ and $\pi_0,\dots ,\pi_k $.
Denote by $\sigma $ this new permutation which keeps track of the order
of our $m$ blocks as they are progressively refined.
We initialize $\sigma $ to the
identity (so blocks are initially ordered left to right).  Each time
a block $B$ is refined, 
$B$ is shifted to a later position in the evolving block order $\sigma $.  
Specifically, 
$B$ is shifted to the last position among blocks which are  
similar to $B$ in the sense of the following interwoven definitions.

\begin{defn}
{\rm
A series of consecutive covering relations $u_0\prec\cdots\prec u_{km}$ is
called a {\bf similarity series} if there is some collection of blocks
$B_1,\dots ,B_m$ that are similar at $u_0$ and have each been split
in the same fashion in the saturated chain from $\hat{0}$ to $u_0$.  We
also require for each $0\le i< m$ that the covering relations 
$u_{ik}\prec\cdots\prec u_{(i+1)k}$ split the block $B_{i+1}$ in a
fashion that avoids topological descents (and also avoid ranks that 
would be included in the minimal face for 
a partitioning restricted to $B_{i+1}$).  Furthermore, the blocks
must be split in identical fashion within a similarity series.  
}
\end{defn}

\medskip
The requirement about avoiding ranks that would be included in a 
partitioning minimal face is discussed more just prior to Theorem
~\ref{biglinks}.  By definition, similarity series' are non-overlapping.
 
\begin{defn}
{\rm
Let us define {\bf similarity} of blocks recursively as follows.  When
$m$ consecutive ranks insert bars from left to right in a single block 
creating $m$ left children of equal size, these children are at this point
all similar.  A collection 
of blocks $B_1,\dots ,B_t$ which are similar at $u$ will still be similar
at $v$ for $u<v$ if every time any one of the blocks $B_i$ 
appears in the interval from $u$ to $v$, 
it appears as part of a similarity series for $B_1,\dots ,B_t$
(though this similarity series might continue beyond $v$ or begin prior to
$u$).
}
\end{defn}

Notice that a block $B_i\in \{ B_1,\dots ,B_m \} $ 
ceases to be similar to the other blocks 
when two bars are inserted in consecutive steps in 
$B_i$ as a topological descent.  Similarity is also broken at ranks
which for some other reason are included in the 
minimal face assigned to the facet restricted to $B_i$ in the 
partitioning for the complex restricted to $B_i$.  

As an example, each facet in Figure ~\ref{j2} has two 
similar blocks at rank 2 which remain similar at rank 6.  This similarity
would have been broken if 
bars were instead inserted from right to left at ranks $3,4$ 
or $5,6$.  Notice that
when two similar blocks are split in identical
fashion, but there are intermediate steps splitting other blocks in 
other ways, then block similarity is broken.  Once $B_i$ and $B_j$ cease
to be similar, $\sigma $ henceforth
preserves the relative order of $B_i$ and $B_j$.

At each refinement step, $\sigma $ orders blocks within each similarity class 
according to the order in which they have most recently been refined.
The descents in $\sigma_{final}^{-1}$ determine which splitters to 
include in the minimal face associated to a facet.  
At any particular rank, $\sigma $ reflects the 
partial evolution based on block refinement up to this point.  
Now let us define what it means for two consecutive refinements
$u\prec v, v\prec w$ to be increasing in the relative transpose order. 

\begin{defn}
Let us represent each of the $km$ positions into which 
bars may be inserted by a pair 
$(B,s)$ consisting of the block $B$ being split and 
the slot $s$ within $B$.  The {\bf relative transpose order} on 
bar positions $(B,s)$ satisfies $(B,s)<(B',s')$ if (1) $s<s'$ or
(2) $s=s'$ and $\pi(B)<\pi(B')$, evaluating $\pi $ at $u$, i.e. 
just prior to both of the labels to be compared.
\end{defn}

Descents in the relative transpose order indicate which ranks to 
include in assigning minimal faces to facets in an analogous fashion
to Section ~\ref{linkpart}.
As an example, the first facet in Figure ~\ref{j2} has a descent
at rank 4 in the relative transpose order, but not at ranks 3 or 5, 
leading us to assign the minimal face of support $1,4,7$ to the first
facet of Figure ~\ref{j2}.

When two bars are inserted in consecutive steps into a block 
$B_i\in\{ B_1,\dots ,B_m\} $, 
we turn to a partitioning for the quotient complex
given by restriction to $B_i$ to decide whether to include the rank in the 
minimal face assigned to our facet.  In particular, the rank
is included as topological descent 
unless the restriction gives a non-shelling step in the lexicographic
order on $\Delta(\Pi_{|B_i|})/S_{|B_i|}$, in which case we turn to 
embedded instances of partitioning link subcomplexes.  

\begin{thm}\label{biglinks}
The link subcomplex $\Delta $ defined in this section is partitionable.
\end{thm}
\proof
The minimal face $G_j$ assigned to the facet $F_j$
is obtained by restricting to the following ranks of $F_j$.
\begin{itemize}
\item
any slot rank separating consecutive bar insertions into distinct blocks 
$B_i,B_j\in\{ B_1,\dots ,B_m\} $ which is a descent in the relative 
transpose order on bar positions.
\item
any slot rank separating consecutive bar insertions into a single
block $B_i\in\{ B_1,\dots ,B_m \} $ which 
would be included in the minimal face assigned to the 
facet restricted to $\Delta(\Pi_{|B_i|})/S_{|B_i|}$ in 
the partitioning of the 
quotient complex $\Delta(\Pi_{|B_i|})/S_{|B_i|}$.
\item
any splitter such that $\sigma_{final}^{-1}$
has a descent at the splitter location.
\end{itemize}

With these choices, the arguments of Section ~\ref{linkpart} generalize.
Namely, each face $F$ extends by increasing chains in
the relative transpose order to a unique facet $F_i$ such
that $F\in [G_i,F_i]$.  Begin by splitting
the leftmost of equivalent blocks at each rank.  Then extend to a
facet $\overline{F}$
that is increasing in the relative transpose order on each
interval, except perhaps for descents from wrap-around.  Next, 
permute the blocks in between consecutive splitters included in $F$ so 
that $\sigma_{final}$ is increasing in between splitters.  Now extend
each interval in the unique way that is increasing in the relative
tranpose order to obtained the desired $F_i$ 
with the property that $F\in [G_i,F_i]$.
\EOP

\begin{figure}[h]
\begin{picture}(180,42)(48,9)
\psfig{figure=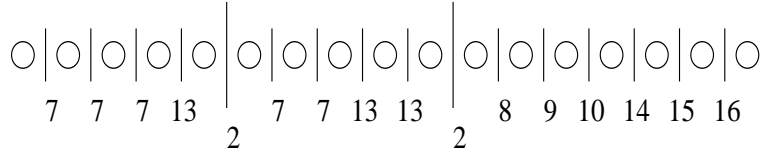,height=2cm,width=10cm}
\end{picture}
\caption{A face with nonconsecutive steps splitting equivalent blocks}
\label{nonconsec}
\end{figure}

\medskip
More generally, we will also need to account for
non-shelling steps in which a set of similar blocks is partially
refined, then other blocks are refined before 
further refining the set of equivalent blocks.  Figure ~\ref{nonconsec}
gives an example of a face from such a non-shelling step.  
Partitioning must accomodate any number of such alternations.

\begin{rk}\label{alternation}
\rm
Partitioning link subcomplexes in which 
the similar blocks need not be refined entirely in consecutive steps
is done as in Theorem ~\ref{biglinks}.  
That is, whenever, two steps refining similar
blocks are separated by a step refining an outside block $B$, there is a 
topological descent either in departing from the collection of 
similar blocks or in returning to it, because the similar blocks occur 
consecutively (as will be verified in Lemma ~\ref{consecblocks}).  This
descent breaks the similarity
of two blocks when one is refined before the refinement of $B$ and 
the other is refined 
afterwards.  This is the only amendment to the argument of 
Theorem ~\ref{biglinks} which is needed.
\end{rk}

\subsection{Extension to a partitioning for $\dps $}

Let us now characterize all non-shelling steps and show how to merge the 
partitionings of Section ~\ref{nextlink} into a partitioning
for all of $\dps $.
By Theorem ~\ref{icc}, a 
non-shelling step may only arise when skipping some minimal
rank set $i_1,\dots ,
i_k$ renders two blocks $B, B'$ in a facet $F_j$ equivalent, and the block
$B$ or $B'$ which is farther to the right is split first by $F_j$.
By the increasing chain condition (cf. Theorem ~\ref{icc}), 
the blocks $B,B'$ must be created within a sequence of consecutive
refinement steps that create a set of identical 
blocks from left to right from a single parent.  

\begin{lem}\label{consecblocks}
Similar blocks are positioned consecutively.
\end{lem}
\proof
Consecutiveness follows from the fact that 
the blocks must be created from a single parent in consecutive steps with 
no topological descents; the one case requiring special 
care is when the rightmost of similar blocks 
is larger than the other block created in the same step, so then 
a smaller block comes in 
between $m-1$ of the similar blocks and the last one, but 
this gives a topological descent immediately before the 
bar splitting off the last of the identical blocks,
rendering this last block not equivalent to the others.
\EOP

\medskip
When the last two blocks created (in a single
step) have equal size, then these two blocks are not just similar, but 
actually equivalent.  This implies that the bar between them cannot be a
splitter, in that the inclusion of this rank does not distinguish the two
ranks.  By convention, we always must refine the left of
these two equivalent blocks first, until their equivalence is broken by 
the inclusion of a slot rank.  When the last slot does not have this
rightmost block filled last, then this gives a descent in the inverse to 
the wraparound permutation.  By convention, we then include in the
associated minimal face the rank in the first slot
(rather than the forbidden rightmost splitter in the 0-slot).

Because the set of $m$ similar blocks appear 
consecutively, any two facets sharing 
a link of the type considered in Section ~\ref{nextlink} 
will have the same topological
descents among the ranks outside the link.  This ensures
that all faces in a particular link subcomplex will include the same
external ranks in their minimal faces both from descents and also
from other external partitionings.
In conclusion, we get a partitioning for $\dps $, as confirmed next.

\begin{thm}
The quotient complex $\Delta (\Pi_n)/S_n$ is partitionable.
\end{thm}
\proof
For each facet $F_j$ which is inserted as a lexicographic shelling step
as in Section 4, 
let $G_j$ be the face comprised of the topological descents in $F_j$.
By Theorem ~\ref{mult-face}, every 
nonshelling step $F_j$ has at least one series of consecutive
bar insertions creating similar blocks by inserting splitters from
left to right.  For each such collection of similar blocks,
use the partitioning of the link subcomplex with the appropriate 
$m,k$ and with the partial refinements that reflect alternation between 
refining the $m$ blocks and the blocks outside the segment (cf. Remark 
~\ref{alternation}) to determine which
ranks to include in $G_j$ among the ranks belonging to the link.

Among the ranks of $F_j$ not belonging to any such link subcomplex, include 
in $G_j$ the topological descents according to the 
lexicographic order on facets given in Section 4.
To be sure this this assignment of minimal faces to facets 
gives a partitioning, we need the fact that 
facets involved in the same link 
the same topological descents outside of the link.  
This follows from the
characterization of non-shelling steps, and in particular from the 
fact that equivalent blocks must occur consecutively.  Thus, any fixed 
bar insertion outside the collection of equivalent blocks will either be 
to the right of all possible bar positions within
the collection of equivalent blocks or to the left of them all.  
Thus, the non-shelling steps
collectively contribute exactly the faces needed to complement the 
lexicographic shelling steps.  This follows because for 
any non-shelling step, all the minimal 
faces that first appear lexicographically in it
are obtained by taking a minimal face in the partitioning of a link
subcomplex, as in Section ~\ref{nextlink}, and adding the ranks outside
the link which are topological descents or are chosen in another link
subcomplex partitioning.
\EOP

\section{Obtaining the multiplicity of the trivial representation 
from a quotient complex shelling or partitioning}\label{trivsec}

Following [St2],
denote by $\alpha_S (\Pi_n )$ the symmetric group action on chains with
rank set $S$.  Let $\beta_S (\Pi_n )$ be the symmetric group
representation on top homology of the rank-selected partition lattice
with ranks belonging to $S$.  Information about the flag $f$-vector and 
flag $h$-vector may be found (for instance) in [St4].  Recall, 
$$ \beta_S  = \sum_{T\subseteq S} (-1)^{|S-T|} \alpha_T $$
where $\beta_S $ is the induced representation on rank-restricted homology,
since $\Pi_n $ is Cohen-Macaulay.
Note that the multiplicity $\langle \alpha_T ,1\rangle $
of the trivial representation in $\alpha_T(\Pi_n )$
is the number of orbits in the induced representation 
on chains with rank set $T$, namely
$\langle \alpha_T ,1\rangle = 
f_T (\Delta (\Pi_n )/S_n )$.  This together with inner product linearity
yields
\begin{eqnarray*}
\langle \beta_S ,1\rangle 
&=& \sum_{T\subseteq S} (-1)^{|S-T|} \langle\alpha_T ,1\rangle \\
&=& \sum_{T\subseteq S} (-1)^{|S-T|} f_S (\Delta (\Pi_n )/S_n ) \\
&=& h_S (\Delta (\Pi_n )/S_n )
\end{eqnarray*}
where $h_S$ is the term indexed by $S$ in the flag $h$-vector for $\dps $.
Recall that 
$h_S (\Delta )$ counts the minimal new faces colored by $S$ in any shelling (or
partitioning) for a balanced boolean cell complex.
Thus, a shelling yields an interpretation for the multiplicity $\langle 
\beta_S ,1 \rangle $
of the trivial representation in $\beta_S (\Pi_n)$.  
Furthermore, the collection of
colors in the 
minimal new face in a lexicographic shelling step is the collection
of ranks where (topological) descents occur in the facet being inserted.

The above discussion applies to any rank-preserving group action on any finite,
ranked Cohen-Macaulay poset with a shelling or partitioning.
Thus, we also get an interpretation for the multiplicity
of the trivial representation in the 
$S_2\wr S_n$ action on the homology of the rank-selected boolean lattice
$B_{2n}$ from our shelling for $\Delta (B_{2n})/S_2\wr S_n$.

\section{Predicting when $b_S(n) =0$}\label{bs0sec}

The partitioning for $\dps $ yields a combinatorial 
interpretation for $b_S $ as the number of 
facets with minimal new face 
colored by $S$ in the partitioning for $\dps $.
The following appeared as Conjecture 4.11 in [BK].  

\begin{conj}[Babson-Kozlov]\label{bk}
The rank-selected quotient complex $\Delta(\Pi_n^S)/S_n$ is contractible
if and only if $b_S(n)=0$.
\end{conj}

Notice that $b_S(n)>0$ if and only if the partitioning for 
$\Delta (\Pi_n)/S_n$ has a minimal face of support $S$.  A shelling for
$\Delta (\Pi_n)/S_n$ would resolve the conjecture affirmatively.

Next, we recover a result of Hanlon [Ha] (see [Su], [Ko1] for proofs by
other methods).  Our argument is included to give a concrete
example of how the partitioning for $\Delta (\Pi_n)/S_n$ may lead to 
results about $\langle \beta_S ,1\rangle $.  More extensive results
of this nature, including the confirmation of two conjecture from [Su], 
have recently been developed in a joint project with Phil
Hanlon [HH].  

Our proof below resembles the one given by Kozlov
in [Ko1] in the sense that we both show that 
$h_{1,\dots ,i} (\Delta (\Pi_n)/S_n) =0$
and use the fact that $h_S(\Delta (\Pi_n)/S_n)=\langle \beta_S ,1\rangle $.
However, Kozlov's proof is topological and ours is combinatorial.
More specifically, Kozlov uses the 
interpretation of $h_S (\Delta )$ as the reduced
Euler characteristic for $\Delta^S$, deducing that it is 0 from the fact that 
$\Delta (\Pi_n^{1,\dots ,i})/S_n$ is collapsible.  We instead
use the fact that $h_S (\Delta )$ counts boolean intervals with minimal 
element of support $S$ in a partitioning for $\Delta $.

\begin{figure}[h]
\begin{picture}(100,55)(20,5)
\psfig{figure=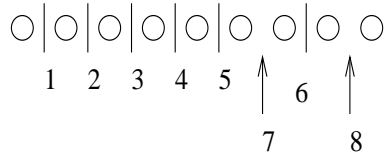,height=2cm,width=5cm}
\end{picture}
\caption{The impossibility of descent set $\{ n-i-1 ,n-i,\dots ,n-1\} $}
\label{hanlon}
\end{figure}

To show $b_S(n)=0$ for a particular $S$, we will show
that $S$ never occurs as a descent set in a shelling step or as the 
support of the minimal new face assigned to some non-shelling step.  
We work in terms of the dual poset to 
the partition lattice as considered in [Ha] and [St2].  To account for this,
we must replace each $i\in S$ by $n-i-1$.  
Thus, the next theorem requires that we show there are no shelling steps
that begin with entirely (topological) ascents and then consist
entirely of (topological) descents as soon as the first (topological) 
descent occurs; we must also check a 
similar condition for the non-shelling steps.
\begin{thm}
The multiplicity $b_S(n)$
of the trivial representation in 
$\beta_S$ is 0 for $S=\{ 1,\dots ,i \} $ and $1\le i\le n-2$.
\end{thm}
\proof
First consider the facets that are inserted as shelling steps.
Recall our convention
of placing bars as far to the left as possible at each 
step.  A facet with descent set $n-i-1,n-i,\dots ,n-3,n-2$
could only come from initially inserting bars left to
right creating blocks nondecreasing in size; after this, we would
need to completely refine the rightmost
block using only topological descents and then proceed similarly 
through the remaining blocks from right
to left.  If the 
initial insertion of bars left to right created any blocks of size 
greater than 2, it
would be impossible to later refine such a block with only 
topological descents.  Thus,
we may assume bars are initially placed left to right creating blocks 
entirely of size
1 followed by blocks entirely of size 2, as in Figure \ref{hanlon}.
However, notice that the two rightmost blocks of size 2 will be equivalent
to each other, which means the left one must be split first.  This 
necessitates the existence of a (topological) ascent at some point 
after the first (topological) descent, as in Figure \ref{hanlon}, so we are
done with the shelling steps.

The non-shelling steps all have support including ranks other than only
the final string by 
virtue of creating equivalent blocks to be partitioned from left to right.
A descent is necessitated by the insertion of a bar creating the 
rightmost of the equivalent blocks farther to the right than bars to be
inserted within the equivalent blocks; there will be an ascent at some point
after this descent since the partitionings of link subcomplexes 
cannot have facets with minimal faces of full support.
\EOP

\section{An application to sub-rings of invariant polynomials}\label{invtsec}

Let us 
review notation from [GS] needed to state a theorem of Garsia and
Stanton (which appears as Theorem 6.2 in [GS] and is restated as Theorem
\ref{gs} below).  Let $C$ be a 
balanced boolean cell complex consisting of vertices 
$x_1,\cdots , x_n$, ordered in a way that is compatible with their 
colors; that is, we choose a vertex order such that
$x_i$ is colored with a smaller number
than $x_j$ for all $1\le i<j \le n$.  If $c$ is a face of $C$ consisting
of vertices $x_{i_1},\cdots ,x_{i_k}$, then 
denote by $x(c)$ the monomial $x_{i_1}\cdots x_{i_k}$.
When a group $H$ acts on $C$, let
$$R^H (x(c)) = \frac{1}{|H|}\sum_{h\in H} hx(c) = \frac{1}{|H|}\sum_{h\in H}
x(hc).$$  Let $\theta_i = \sum_{c(v)=i} x(v)$, a sum over vertices in $C$ 
of color $i$.
A set of chain monomials $\{ x(b) | b\in B\} $ given by a collection $B$ 
of chains in a poset $P$ is called a basic set if every
element $Q$ of the Stanley-Reisner ring $R_P$ has a unique expression
$$Q = \sum_{b\in B} x(b) Q_b (\theta_1,\dots ,\theta_d)$$ where the 
coefficients $Q_b(\theta_1,\dots ,\theta_d )$ are polynomials with rational
coefficients in the variables $\theta_1,\cdots ,\theta_d $.  
All Cohen-Macaulay posets have such
basic sets, and [GS] shows how a quotient complex
shelling gives an explicit 
basic set.  

\begin{thm}[Garsia and Stanton]\label{gs}
If $C/H$ has a shelling $F_1,\dots ,F_k $ where $R(F_i) $ is the minimal
new face in $F_i$ and $b_i$ is a representative within $C$ 
of the orbit of $R(F_i)$,
then the orbit polynomials $R^H x(b_i) $ form a basic set for
$R^H R_P$.  
\end{thm}

Thus, our lexicographic shelling for 
$\dwreath $ yields a ring invariant basic set.
A simple description of which descent sets
occur in the lexicographic shelling would be desirable since it would
yield a nice description of 
these ring invariant basic sets, according to Theorem \ref{gs}.  It 
would also be interesting to determine whether the incidence matrix
given by our partitioning for $\dps $ is nonsingular, since that would
imply a basic set for $k[\Delta(\Pi_n)]^{S_n}$, by another result of [GS].

Stanley showed in [St1] that the face ring of a 
Cohen-Macaulay simplicial poset is a Cohen-Macaulay ring.  According
to Proposition \ref{goofcm}, our shelling
for $\Delta(B_{2n})/S_2\wr S_n$ implies that
the face ring $k[\Delta(B_{2n})/S_2\wr S_n ]$ 
is Cohen-Macaulay (over the integers).
However, Reiner constructed an isomorphism in [Re] between the face ring
$k[\Delta(P)/G]$ and the ring of invariants $k[\Delta(P)]^G$. 
Reiner's result allows us to conclude that the invariant subring
$k[B_{2n}]^{S_2\wr S_n}$ is Cohen-Macaulay over the integers.

Furthermore, Reiner showed that in the case of quotients of the 
Boolean algebra by a permutation subgroup $G$ of $S_n$, 
the invariant subring of
the polynomial ring $k[x_1,\dots ,x_n]$ under the action of $G$
is Cohen-Macaulay 
over $R$ whenever the same is true for the 
invariant subring $R[B_n]^G$
of the Boolean algebra's face ring over $R$; Reiner recently provided
a proof for this (formerly unpublished) result as an appendix to [He].
Since $S_2\wr S_n$ is such a permutation subgroup of $S_{2n}$,
we may conclude the following from our shelling for $\Delta(B_{2n})/
S_2\wr S_n$. 

\begin{thm}
The subring of invariant polynomials 
$k[x_1,\dots ,x_{2n}]^{S_2\wr S_n }$ is 
Cohen-Macaulay for any field $k$.
\end{thm}

When $k$ is the field of complex numbers, this is a special case
of a result from [HE], but the shelling gives the 
Cohen-Macaulay property also for fields of finite characteristic, or
equivalently for integer coefficients. 

\section*{Acknowledgments}

The author thanks Vic Reiner for suggesting
$\Delta(\Pi_n)/S_n$ and $\dwreath $ as candidates for the 
lexicographic shelling criterion for balanced complexes as well
as for explaining to her the potential implications discussed 
in Sections 6 and 8 and also suggesting the use of 
Molien's Theorem in conjunction with Section 3.  She is also grateful 
to Eric Babson, Phil Hanlon, Robert Kleinberg and Richard Stanley
for helpful discussions, Axel Hultman for proofreading portions of
this paper and the anonymous referees for their generous advice on how
to improve various aspects of this paper.


\begin{thebibliography}{23}

\bibitem[BK]{BK}
E. Babson and D. Kozlov, {\it Group actions on posets}.  Preprint 1998.

\bibitem[Bj1]{B1}
A. Bj\"orner, {\it Posets, regular CW complexes and Bruhat order}, European
J. Combin. {\bf 5, no. 1} (1984), 7-16

\bibitem[Bj2]{B}
A. Bj\"orner, {\it Shellable and Cohen-Macaulay partially ordered sets},
Trans. Amer. Math. Soc. {\bf 260, No. 1} (1980), 159-183.

\bibitem[Bj3]{Bj3}
A. Bj\"orner,  Topological Methods, in
{\it Handbook of Combinatorics} (R. Graham, M. Gr\"otschel and L Lovasz, eds.),
North-Holland, Amsterdam, 1993.

\bibitem[BW1]{BW1}
A. Bj\"orner and M. Wachs, {\it On lexicographically shellable posets},
Trans. Amer. Math. Soc. {\bf 277, No. 1} (1983), 323-341.

\bibitem[BW2]{BW2}
A. Bj\"orner and M. Wachs, {\it Bruhat order of Coxeter groups and shellability},
Adv. in Math. {\bf 43} (1982), 87-100.

\bibitem[Du]{Du}
A. Duval, {\it Free resolutions of simplicial posets}, Journ. of Algebra {\bf 188} (1997),
363-399.

\bibitem[GS]{GS}
A. Garsia and D. Stanton, {\it Group actions of Stanley-Reisner rings and 
invariants of permutation groups}, Adv. in Math. {\bf 51} (1984), no. 2, 
107-201

\bibitem[Ha]{Ha}
P. Hanlon, {\it A proof of a conjecture of Stanley concerning partitions of a set}, 
European J. Combin. {\bf 4, no. 2} (1983), 137-141.

\bibitem[HH]{HH}
P. Hanlon and P. Hersh, {\it Multiplicity of the trivial representation in 
rank-selected homology of the partition lattice}, To appear in J. Algebra.

\bibitem[He]{He}
P. Hersh, {\it A partitioning and related properties for the quotient complex
$\Delta(B_{lm})/S_l\wr S_m$}, To appear in J. Pure and Appl. Alg.

\bibitem[HK]{HK}
P. Hersh and R. Kleinberg, {\it The refinement complex of the poset of 
partitions of a multiset}, In preparation.

\bibitem[HE]{HE}
M. Hochster and J.A. Eagon, {\it Cohen-Macaulay rings, invariant theory, and 
the generic perfection of determinantal loci}, Amer. J. Math. {\bf 93} (1971), 1020-1058.

\bibitem[Hu]{Hu}
A. Hultman, {\it Lexicographic shellability and quotient complexes}, 
To appear in J. Algebraic Combinatorics.

\bibitem[Ko1]{Ko1}
D. Kozlov, {\it Collapsibility of $\Delta(\Pi_n)/S_n $ and some related 
CW complexes}, Proc. Amer. Math. Soc. {\bf 128}, (2000), no. 8, 2253-2259.

\bibitem[Ko2]{Ko2}
D. Kozlov,  {\it General lexicographic shellability and orbit 
arrangements}, Ann. ~of Comb.~{\bf 1 (1)}, (1997), 67-90.

\bibitem[Ko3]{Ko3}
D. Kozlov, {\it Rational homology of spaces of complex monic polynomials
with multiple roots}, math CO/0111167.

\bibitem[Mu]{Mu}
J. Munkres, {\it Topological results in combinatorics}, Michigan
Math. J. {\bf 31} 113-128.

\bibitem[Re]{Re}
V. Reiner, {\it Quotients of Coxeter complexes and $P$-partitions},
Memoirs Amer. Math. Soc. {\bf 95}, January 1992.

\bibitem[St1]{St1}
R. Stanley, {\it Invariants of finite groups and their applications to 
combinatorics}, Bull. Amer. Math. Soc. {\bf 1} (1979), 475-511.

\bibitem[St2]{St2}
R. Stanley, {\it Some aspects of groups acting on finite posets}, J. Combin.
Theory Ser. A, {\bf 32 (2)} (1982), 132-161.

\bibitem[St3]{St3}
R. Stanley, {\it $f$-vectors and $h$-vectors of simplicial posets}, J. Pure
and Applied Algebra, {\bf 71} (1991), 319-331.

\bibitem[St4]{St4}
R. Stanley, Combinatorics and Commutative Algebra, second ed., Birkh\"auser,
Boston, 1996.

\bibitem[Su]{Su}
S. Sundaram, The homology representations of the symmetric group on 
Cohen-Macaulay subposets of the partition lattice, {\it Adv. Math.}
{\bf 104} (1994), 225-296.

\bibitem[Wa]{Wa}
M. Wachs, {\it A basis for the homology of the $d$-divisible 
partition lattices}, Adv. in Math. {\bf 117}, no. 2 (1996) 294-318.
 
\bibitem[Zi]{Zi}
G. Ziegler, {\it On the poset of partitions of an integer}, J. Combin.
Theory Ser. A {\bf 42, No. 2}, (1986) 215-222.
\end{thebibliography}
\end{document}